\documentclass[10pt,letterpaper,twoside,english,aps,manuscript,nofootinbib, nobibnotes, secnumarabic, notitlepage]{revtex4}
\usepackage[T1]{fontenc}
\usepackage[latin9]{inputenc}
\usepackage{xcolor}
\usepackage{babel}
\usepackage{units}
\usepackage{amsthm}
\usepackage{amsmath}
\usepackage{amssymb}
\usepackage{graphicx}
\PassOptionsToPackage{normalem}{ulem}
\usepackage{ulem}
\usepackage[unicode=true]
 {hyperref}
\usepackage{breakurl}

\makeatletter

\providecommand{\tabularnewline}{\\}
\newcommand{\lyxdot}{.}

\providecolor{lyxadded}{rgb}{0,0,1}
\providecolor{lyxdeleted}{rgb}{1,0,0}

\@ifundefined{textcolor}{}
{%
 \definecolor{BLACK}{gray}{0}
 \definecolor{WHITE}{gray}{1}
 \definecolor{RED}{rgb}{1,0,0}
 \definecolor{GREEN}{rgb}{0,1,0}
 \definecolor{BLUE}{rgb}{0,0,1}
 \definecolor{CYAN}{cmyk}{1,0,0,0}
 \definecolor{MAGENTA}{cmyk}{0,1,0,0}
 \definecolor{YELLOW}{cmyk}{0,0,1,0}
}
\theoremstyle{plain}
\newtheorem{thm}{\protect\theoremname}[section]
  \theoremstyle{remark}
  \newtheorem{rem}[thm]{\protect\remarkname}

\usepackage{multirow}
\usepackage{appendix}

\usepackage{placeins}

\def\thetable{\@arabic\c@table}

\def\frontmatter@abstractheading{}

 \numberwithin{equation}{section}

\AtBeginDocument{
  
}

\makeatother

  \providecommand{\remarkname}{Remark}
\providecommand{\theoremname}{Theorem}

\begin{document}

\title{All-Possible-Couplings Approach \\to \\Measuring Probabilistic
Context}

\author{Ehtibar N. Dzhafarov\textsuperscript{1*} and Janne V. Kujala\textsuperscript{2}}

\affiliation{$\,$}

\affiliation{\textsuperscript{1}Department of Psychological Sciences, Purdue
University, USA,\\}

\affiliation{\textsuperscript{2}Department of Mathematical Information Technology,
University of Jyväskylä,\\\\ Finland\\\textsuperscript{*}Corresponding
author (email: ehtibar@purdue.edu)}
\begin{abstract}
\markboth{Dzhafarov and Kujala}{Measuring Context}
\end{abstract}
\maketitle
$\,$\vspace{0.05\textheight}

\begin{center}
\textbf{\large ABSTRACT}
\par\end{center}{\large \par}

\noindent From behavioral sciences to biology to quantum mechanics,
one encounters situations where (i) a system outputs several random
variables in response to several inputs, (ii) for each of these responses
only some of the inputs may \textquotedblleft{}directly\textquotedblright{}
influence them, but (iii) other inputs provide a \textquotedblleft{}context\textquotedblright{}
for this response by influencing its probabilistic relations to other
responses. These contextual influences are very different, say, in
classical kinetic theory and in the entanglement paradigm of quantum
mechanics, which are traditionally interpreted as representing different
forms of physical determinism. One can mathematically construct systems
with other types of contextuality, whether or not empirically realizable:
those that form special cases of the classical type, those that fall
between the classical and quantum ones, and those that violate the
quantum type. We show how one can quantify and classify all logically
possible contextual influences by studying various sets of probabilistic
couplings, i.e., sets of joint distributions imposed on random outputs
recorded at different (mutually incompatible) values of inputs. 

\medskip{}


\pagebreak{}

\section{Introduction}

Consider a system with two inputs, $\alpha,\beta$, and two random
outputs, $A,B$, about which it is assumed that $A$ is \emph{not}
influenced by $\beta$, nor $B$ by $\alpha$.\textcolor{black}{{} }A
necessary condition for this selectivity of influences is \emph{marginal
selectivity} \textcolor{black}{{[}}\textcolor{black}{\emph{\ref{enu:Townsend,-J.T.,-&}}}\textcolor{black}{{]}}:
changes in the values of $\beta$ do not influence the distribution
of $A$, and analogously for $\alpha$ and $B$. Let, for example,
both inputs and outputs be binary: $\alpha=\left\{ \alpha_{1},\alpha_{2}\right\} $,
$\beta=$$\left\{ \beta_{1},\beta_{2}\right\} $, and $A,B$ attain
values $+1$ and $-1$ each. Denoting by $A_{ij}$ and $B_{ij}$ the
two outputs conditioned on $\alpha=\alpha_{i},\beta=\beta_{j}$ ($i,j\in\left\{ 1,2\right\} $),
the distribution of $\left(A_{ij},B_{ij}\right)$ is described by
the joint probabilities $p_{ij},q_{ij},r_{ij},s_{ij}$ (summing to
1) in the matrix 
\begin{equation}
\begin{array}{|c|c|c|}
\hline \alpha_{i},\beta_{j} & B_{ij}=+1 & B_{ij}=-1\\
\hline A_{ij}=+1 & p_{ij} & q_{ij}\\
\hline A_{ij}=-1 & r_{ij} & s_{ij}
\\\hline \end{array}.\label{eq:matrix}
\end{equation}
Assuming all four combinations $\left\{ \alpha_{1},\alpha_{2}\right\} \times\left\{ \beta_{1},\beta_{2}\right\} $
are possible, marginal selectivity in this example means 
\begin{equation}
\begin{array}{c}
p_{i1}+q_{i1}=p_{i2}+q_{i2}=\Pr\left[A_{ij}=+1\right],\\
p_{1j}+r_{1j}=p_{2j}+r_{2j}=\Pr\left[B_{ij}=+1\right],
\end{array}\label{eq:marginal selectivity}
\end{equation}
for all $i,j\in\left\{ 1,2\right\} $. 

The assumption of \emph{selective influences}, however, is stronger.
It requires that the joint distribution of the two outputs satisfies,
for all $i,j\in\left\{ 1,2\right\} $, 
\begin{equation}
\left(A_{ij},B_{ij}\right)\sim\left(f\left(R,\alpha_{i}\right),g\left(R,\beta_{j}\right)\right)\label{eq:selective influences}
\end{equation}
where $\sim$ stands for ``has the same distribution as,'' $f,g$
are some functions, and $R$ is a source of randomness that does not
depend on $\alpha,\beta$ {[}\emph{\ref{enu:Bell,-J.-(1964).}-\ref{enu:Dzhafarov,-E.N.-(2003).}}{]}.
In our example (\ref{eq:matrix}) this means
\begin{equation}
p_{ij}=\Pr\left[f\left(R,\alpha_{i}\right)=+1,g\left(R,\beta_{j}\right)=+1\right],\; r_{ij}=\Pr\left[f\left(R,\alpha_{i}\right)=+1,g\left(R,\beta_{j}\right)=-1\right],\;\textnormal{etc.}
\end{equation}
In the quantum mechanical context (see below) $R$ is interpreted
as ``hidden variables.'' Such a representation may or may not exist
when marginal selectivity is satisfied. For instance, the latter is
satisfied in the following four distributions,
\begin{equation}
\begin{array}{ccc}
\begin{array}{|c|c|c|}
\hline \alpha_{1},\beta_{1} & B_{11}=+1 & B_{11}=-1\\
\hline A_{11}=+1 & \nicefrac{1}{4} & 0\\
\hline A_{11}=-1 & 0 & \nicefrac{3}{4}
\\\hline \end{array} &  & \begin{array}{|c|c|c|}
\hline \alpha_{1},\beta_{2} & B_{12}=+1 & B_{12}=-1\\
\hline A_{12}=+1 & 0 & \nicefrac{1}{4}\\
\hline A_{12}=-1 & \nicefrac{1}{2} & \nicefrac{1}{4}
\\\hline \end{array}\\
\\
\begin{array}{|c|c|c|}
\hline \alpha_{2},\beta_{1} & B_{21}=+1 & B_{21}=-1\\
\hline A_{21}=+1 & 0 & \nicefrac{1}{2}\\
\hline A_{21}=-1 & \nicefrac{1}{4} & \nicefrac{1}{4}
\\\hline \end{array} &  & \begin{array}{|c|c|c|}
\hline \alpha_{2},\beta_{2} & B_{22}=+1 & B_{22}=-1\\
\hline A_{22}=+1 & 0 & \nicefrac{1}{2}\\
\hline A_{22}=-1 & \nicefrac{1}{2} & 0
\\\hline \end{array}
\end{array}\label{eq:caseNo}
\end{equation}
It can be shown, however, that no representation (\ref{eq:selective influences})
here is possible as the joint probabilities violate the Bell/CHSH
inequalities considered below (Section \ref{sec:Forms-of-context}
and Appendix \ref{sub:Derivation-of-theFine}). At the same time,
a representation in the form of (\ref{eq:selective influences}) is
possible for the similar distributions
\begin{equation}
\begin{array}{ccc}
\begin{array}{|c|c|c|}
\hline \alpha_{1},\beta_{1} & B_{11}=+1 & B_{11}=-1\\
\hline A_{11}=+1 & \nicefrac{1}{4} & 0\\
\hline A_{11}=-1 & 0 & \nicefrac{3}{4}
\\\hline \end{array} &  & \begin{array}{|c|c|c|}
\hline \alpha_{1},\beta_{2} & B_{12}=+1 & B_{12}=-1\\
\hline A_{12}=+1 & \nicefrac{1}{4} & 0\\
\hline A_{12}=-1 & \nicefrac{1}{4} & \nicefrac{1}{2}
\\\hline \end{array}\\
\\
\begin{array}{|c|c|c|}
\hline \alpha_{2},\beta_{1} & B_{21}=+1 & B_{21}=-1\\
\hline A_{21}=+1 & 0 & \nicefrac{1}{2}\\
\hline A_{21}=-1 & \nicefrac{1}{4} & \nicefrac{1}{4}
\\\hline \end{array} &  & \begin{array}{|c|c|c|}
\hline \alpha_{2},\beta_{2} & B_{22}=+1 & B_{22}=-1\\
\hline A_{22}=+1 & 0 & \nicefrac{1}{2}\\
\hline A_{22}=-1 & \nicefrac{1}{2} & 0
\\\hline \end{array}
\end{array}\label{eq:caseYes}
\end{equation}
One can think of $\alpha$ and $\beta$ in (\ref{eq:caseNo}) and
(\ref{eq:caseYes}) as being involved in\emph{ }different kinds of\emph{
probabilistic context} for the ``direct'' dependence of, respectively,
$B$ on $\beta$ and $A$ on $\alpha$. 

We propose a principled way of quantifying and classifying conceivable
contextual influences, whether within or outside the scope of (\ref{eq:selective influences}).
Our approach is neutral with respect to such issues as causality or
what distinguishes direct influences from contextual. We merely accept
as a given a diagram of direct input-output correspondences (e.g.,
$A\leftarrow\alpha,B\leftarrow\beta$) and study the joint distribution
of the outputs at all possible values of the inputs. The interpretation
of the diagram is irrelevant insofar as it is compatible with the
observed pattern of marginal selectivity:as $\alpha$ changes while
$\beta$ remains fixed, the distribution of $B$ does not change,
and as $\beta$ changes while $\alpha$ remains fixed, the distribution
of $A$ does not change. Note that the distribution of $A$ may but
does not have to change in response to changes in $\alpha$, and analogously
for $B$ and $\beta$.

Our approach is maximally general in the sense of applying to arbitrary
sets of inputs and outputs (see Appendix A). To demonstrate it by
detailed computations, however, we confine ourselves in the main text
to binary $\alpha,\beta$ influencing binary $A,B$; and even more
narrowly, to the ``homogeneous'' case with the two values of both
$A$ and $B$ equiprobable at all values of the inputs $\alpha_{i},\beta_{j}$
($i,j\in\left\{ 1,2\right\} $),
\begin{equation}
\Pr\left[A_{ij}=+1\right]=\Pr\left[B_{ij}=+1\right]=\nicefrac{1}{2}.\label{eq:narrow}
\end{equation}
 Marginal selectivity then is satisfied trivially (because all marginal
distributions are fixed).

\begin{figure}
\begin{centering}
\includegraphics[bb=150bp 200bp 592bp 612bp,scale=0.6]{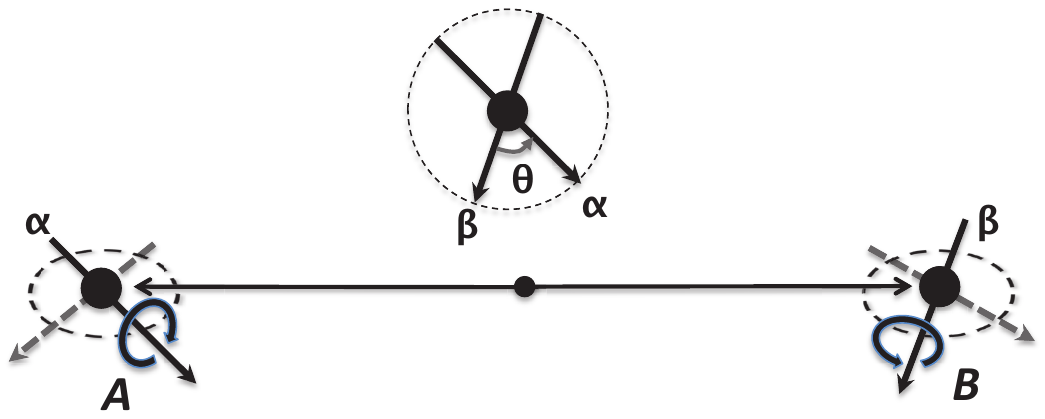}
\par\end{centering}

\caption{Schematic representation of two spin-$\nicefrac{1}{2}$ particles,
e.g., electrons, in the singlet state (represented by $\left\langle \uparrow\right|\otimes\left\langle \downarrow\right|-\left|\downarrow\right\rangle \otimes\left|\uparrow\right\rangle $
in quantum-mechanical notation) running away from each other. The
directions $\alpha$ and $\beta$ are detector settings for spin measurements
(in our language, inputs). The measured spins $A$ and $B$ (outputs)
in these directions are shown by rotation arrows: one direction of
rotation (say, clockwise) represents ``spin-up''$=+1$ in one particle
and ``spin-down''$=-1$ in the other. By the quantum theory, for
any $\alpha,\beta$, $\Pr\left[A=+1,B=+1\right]=\nicefrac{1}{2}\cos^{2}\nicefrac{\theta}{2}$
(equivalently, expected value of $AB$ is $\cos\theta$). The two
measurements are made simultaneously (in some inertial frame of reference).}

\end{figure}

The example focal for this paper is Bohm's version of the Einstein-Podolsky-Rosen
paradigm (EPR/B) {[}\emph{\ref{enu:D.-Bohm,-&}}{]}: a quantum mechanical
system consisting (in the simplest case) of two entangled spin$\textnormal{-}\nicefrac{1}{2}$
particles separated by a space-like interval (see Fig. 1). The two
inputs here are spin measurements on these particles: input $\alpha$
has two values corresponding to spin axes $\alpha_{1},\alpha_{2}$
chosen for one particle, and input $\beta$ has two values corresponding
to spin axes $\beta_{1},\beta_{2}$ for another particle. The two
outputs are spin values recorded: having chosen axes $\alpha_{i}$
and $\beta_{j}$, $i,j\in\left\{ 1,2\right\} $, one records $A_{ij}$
for the first particle and $B_{ij}$ for the second, each being a
random variable with values $+1$ and $-1$.\textcolor{black}{{} (Note
that }the spins of a given particle along two different axes are \emph{noncommuting}
(see Appendix \ref{sub:Derivation-of-theCirelson}), because of which
if one spin value is determined precisely, +1 or -1, the other one
has a nonzero uncertainty. This means that\textcolor{black}{{} }$\alpha_{1},\alpha_{2}$
considered as measurements yielding precise values of spins are mutually
exclusive, and this is the reason $\alpha_{1},\alpha_{2}$ can be
viewed as values of a single input $\alpha$; and analogously for
$\beta_{1},\beta_{2}$ {[}\emph{\ref{enu:E.N.-Dzhafarov,-&2012}},\emph{\ref{enu:E.N.-Dzhafarov,-in press 2}}{]}.\textcolor{black}{)}
Marginal selectivity (\ref{eq:marginal selectivity}) in this context
is known under a variety of other names, such as \textquotedblleft{}parameter
independence\textquotedblright{} and \textquotedblleft{}physical locality\textquotedblright{}
\textcolor{black}{{[}}\textcolor{black}{\emph{\ref{enu:J.-Cereceda,-Quantum}}}\textcolor{black}{{]}}.\textcolor{black}{{}
We confine ourselves to the case (\ref{eq:narrow}), with the two
spin values +1 and -1 being }equiprobable\textcolor{black}{{} for both
}$A_{ij}$\textcolor{black}{{} and }$B_{ij}$\textcolor{black}{. }

\textcolor{black}{Formally equivalent situations are abundant in behavioral
and social sciences {[}}\emph{\ref{enu:Dzhafarov,-E.N.-(2003).}}\textcolor{black}{\emph{,\ref{enu:Kujala,-J.-V.,2008}-\ref{enu:Schweickert,-R.,-Fisher,}}}\textcolor{black}{{]},
where the issue of selective influences was initially introduced in
{[}}\textcolor{black}{\emph{\ref{enu:S.-Sternberg,-The},\ref{enu:J.-T.-Townsend,1984}}}\textcolor{black}{{]},
in the context of information processing architectures. An example
of a system here (from our laboratory) can be a human observer who
adjusts a visual stimulus until it matches in appearance another,
``target'' visual stimulus. Let the latter be characterized by two
properties, $\alpha$ and $\beta$ (e.g., amplitudes of two Fourier-components),
each varying on two levels, $\alpha_{1},\alpha_{2}$ and $\beta_{1},\beta_{2}$.
Denoting by $S_{ij}^{1}$ and $S_{ij}^{2}$ the corresponding properties
(amplitudes) of the adjusted stimulus in response to $\alpha_{i},\beta_{j}$,
we define a binary random output $A_{ij}$ as having the value ``high''$=+1$
or ``low''$=-1$ according as the variable $S_{ij}^{1}$ is above
or below the median of its distribution; output $B_{ij}$ is defined
from $S_{ij}^{2}$ analogously. Marginal selectivity in the form (\ref{eq:narrow})
is ensured here by construction. }

\textcolor{black}{In an example from a biological domain $S_{ij}^{1}$
and $S_{ij}^{2}$ could be activity levels of two neurons tuned to
two stimulus properties, $\alpha$ and $\beta$, respectively. Making
$\alpha$ and $\beta$ vary on two levels each and defining $A_{ij},B_{ij}$
with respect to the medians of $S_{ij}^{1},S_{ij}^{2}$ by the same
rule as above, we get precisely the same mathematical formulation. }

The formal equivalence of these three examples should by no means
be interpreted as a hint at their physical affinity. Unlike in the
EPR/Bohm paradigm, no physical laws prohibit the activity level $A$
of a neuron tuned to stimulus property $\alpha$ from being affected
by stimulus property $\beta$. Similarly, the amplitude $A$ of the
first Fourier component of the adjusted stimulus in the second example
may very well be affected by the amplitude $\beta$ of the second
Fourier component of the target stimulus. Our only claim is that if
these ``secondary'' influences do not change the marginal distributions
of $A$ and $B$ (which in the two examples in question is ensured
by the definition of $A$ and $B$), they can be viewed within the
framework of a formal treatment that also includes the (physically
very different) case of entangled particles.

\section{\label{sec:Forms-of-context}Forms of context (determinism)}

In the following, symbols $i,j,k$ (possibly with primes) always take
on values $1,2$ each, and each of the outputs $A_{ij},B_{ij}$ takes
on values $+1,-1$ with equal probabilities. Representation (\ref{eq:selective influences})
is equivalent to the existence of a jointly distributed system 
\begin{equation}
H=\left(H_{1}^{1},H_{2}^{1},H_{1}^{2},H_{2}^{2}\right),\label{eq:JDC1}
\end{equation}
such that every output pair $A_{ij},B_{ij}$ is distributed as $H_{i}^{1},H_{j}^{2}$;
in symbols, 
\begin{equation}
\left(H_{i}^{1},H_{j}^{2}\right)\sim\left(A_{ij},B_{ij}\right).\label{eq:JDC2}
\end{equation}
As this entails 
\[
H_{i}^{1}\sim A_{ij},\; H_{j}^{2}\sim B_{ij},
\]
all components of $H$ are random variables with equiprobable +1/-1,
and (\ref{eq:JDC2}) reduces to
\begin{equation}
\Pr\left[A_{ij}=+1,B_{ij}=+1\right]=\Pr\left[H_{i}^{1}=+1,H_{j}^{2}=+1\right].
\end{equation}
The existence of $H$ in (\ref{eq:JDC1}) satisfying (\ref{eq:JDC2})
is known as (a special case of) the \emph{Joint Distribution Criterion}
(JDC) {[}\emph{\ref{enu:Fine,-A.-(1982a).},\ref{enu:Fine,-A.-(1982b).},\ref{enu:Dzhafarov,-E.N.,-&2010},\ref{enu:P.-Suppes,-M.1981},\ref{enu:P.-Suppes,-Representation}}{]}.
It follows from (\ref{eq:selective influences}) by 
\begin{equation}
H_{i}^{1}=f\left(R,\alpha_{i}\right),H_{j}^{2}=g\left(R,\beta_{j}\right).
\end{equation}
Conversely, if (\ref{eq:JDC2}) holds for some $H$, then one can
put $R=H$ and 
\begin{equation}
f\left(H,\alpha_{i}\right)=\mathrm{Proj}_{i}\left(H_{1}^{1},H_{2}^{1},H_{1}^{2},H_{2}^{2}\right),\; g\left(H,\beta_{j}\right)=\mathrm{Proj}_{2+j}\left(H_{1}^{1},H_{2}^{1},H_{1}^{2},H_{2}^{2}\right),
\end{equation}
where $\mathrm{Proj}_{k}$ stands for the ``$k$th member'' (in
the list of arguments). The JDC is a deep criterion that provides
a probabilistic foundation for our understanding of the classical
(non)contextuality (or classical determinism in physics). In particular,
it immediately follows from the JDC that if representation (\ref{eq:selective influences})
for $\left(A_{ij},B_{ij}\right)$ exists, the ``hidden variables''
$R$ can always be reduced to a single discrete random variable with
$2^{4}$ possible values (corresponding to the possible values of
$H$).

Using the same notation as above, 
\begin{equation}
p_{ij}=\Pr\left[A_{ij}=+1,B_{ij}=+1\right],
\end{equation}
the JDC in our case (two binary inputs and two binary outputs with
equiprobable values) is equivalent to four double-inequalities
\begin{equation}
0\leq p_{ij}+p_{ij'}+p_{i'j'}-p_{i'j}\leq1\label{eq:Fine}
\end{equation}
with $i\not=i'$, $j\not=j'$ {[}\emph{\ref{enu:Fine,-A.-(1982a).},\ref{enu:Fine,-A.-(1982b).}}{]}.
(See Appendix \ref{sub:Derivation-of-theFine} for a derivation.)
They are often referred to as \emph{the Bell/CHSH inequalities} (in
the homogeneous form), CHSH acronymizing the authors of {[}\emph{\ref{enu:Clauser,-J.-F.,1969}}{]},
although the first appearance of these inequalities dates to {[}\emph{\ref{enu:Clauser,-J.F.-and}}{]}. 

The theory of the EPR/B paradigm predicts and experimental data confirm
violations of the Bell/CHSH inequalities {[}\emph{\ref{enu:A.-Aspect,-P.1981},\ref{enu:A.-Aspect,-P.1982}}{]},
but quantum mechanics imposes its own constraint on the same linear
combinations of probabilities : 
\begin{equation}
\frac{1-\sqrt{2}}{2}\leq p_{ij}+p_{ij'}+p_{i'j'}-p_{i'j}\leq\frac{1+\sqrt{2}}{2}.\label{eq:Cirelson}
\end{equation}
This constraint is known as the Cirel'son inequalities {[}\ref{enu:Cirel'son-BS-(1980)},
\ref{enu:Landau-LJ-(1987)}{]} (see Appendix \ref{sub:Derivation-of-theCirelson}
for a derivation). Since the class of vectors $\left(p_{11},p_{12},p_{21},p_{22}\right)$
that satisfy these double-inequalities include those allowed by (\ref{eq:Fine})
as a proper subset, it is natural to expect that (\ref{eq:Cirelson})
represents some relaxation, or generalization of the JDC. No such
generalization, however, has been previously proposed. Developing
one is the main goal of this paper.

This generalization is not confined to quantum mechanical systems.
In other (e.g., behavioral) applications, one cannot exclude a priori
the possibility of the bounds $m$ and $M$ in 
\begin{equation}
m\leq p_{ij}+p_{ij'}+p_{i'j'}-p_{i'j}\leq M
\end{equation}
being wider than in (\ref{eq:Cirelson}), or falling between the bounds
in (\ref{eq:Fine}) and (\ref{eq:Cirelson}), or being more narrow
than in (\ref{eq:Fine}). One can think of all kinds of other constraints
imposed on the possible values of $\left(p_{11},p_{12},p_{21},p_{22}\right)$,
from confining this vector to one specific value to allowing it to
vary freely. The latter (``complete chaos'') is represented by the
``no-constraint'' constraint
\begin{equation}
-\nicefrac{1}{2}\leq p_{ij}+p_{ij'}+p_{i'j'}-p_{i'j}\leq\nicefrac{3}{2}\label{eq:chaos}
\end{equation}
with $m=-\nicefrac{1}{2}$ attained if one of $p_{11},p_{12},p_{21},p_{22}$
is $\nicefrac{1}{2}$ and the rest are zero, and $M=\nicefrac{3}{2}$
attained if three of $p_{11},p_{12},p_{21},p_{22}$ are $\nicefrac{1}{2}$
and the remaining one is zero. Recall that we only consider the outputs
with equiprobable outcomes, so
\begin{equation}
0\leq p_{ij}\leq\nicefrac{1}{2}.
\end{equation}
All these conceivable constraints on the possible values of $\left(p_{11},p_{12},p_{21},p_{22}\right)$
represent different forms and degrees of contextual influences. It
would be unsatisfactory if all these possibilities, whether or not
empirically realizable, could not be treated within a unified probabilistic
framework including JDC as a special case. We construct such a framework,
based on the classical (Kolmogorov's) theory of probability and the
probabilistic coupling theory {[}\emph{\ref{enu:H.-Thorisson,-Coupling,}}{]}.

\section{Connections}

It is easy to see that for any vector of probabilities $p=\left(p_{11},p_{12},p_{21},p_{22}\right)$
one can find a jointly distributed system of +1/-1 variables 
\begin{equation}
H=\left(H_{11}^{1},H_{11}^{2},H_{12}^{1},H_{12}^{2},H_{21}^{1},H_{21}^{2},H_{22}^{1},H_{22}^{2}\right)\label{eq:extended H}
\end{equation}
such that 
\begin{equation}
\left[\begin{array}{c}
\left(H_{ij}^{1},H_{ij}^{2}\right)\sim\left(A_{ij},B_{ij}\right)\\
\textnormal{i.e.}\\
\Pr\left[H_{ij}^{1}=+1,H_{ij}^{2}=+1\right]=p_{ij}
\end{array}\right]\label{eq:empirical joints}
\end{equation}
for all $i,j$. The JDC then amounts to additionally assuming that
among all such vectors $H$ there is one with 
\begin{equation}
\Pr\left[H_{i1}^{1}\neq H_{i2}^{1}\right]=0,\Pr\left[H_{1j}^{2}\neq H_{2j}^{2}\right]=0,\label{eq:JDCassumption}
\end{equation}
and this is the assumption that is rejected by quantum theory in the
EPR/B paradigm. Once (\ref{eq:JDCassumption}) is explicitly formulated,
however, it becomes clear that it is not the only way of thinking
of $H$. Since $A_{i1}$ and $A_{i2}$ occur under mutually exclusive
conditions, one cannot identify the distribution of $\left(H_{i1}^{1},H_{i2}^{1}\right)$
with that of $\left(A_{i1},A_{i2}\right)$. The latter does not exist
as a pair of jointly distributed random variables. There is therefore
no privileged pairing scheme for realizations of $H_{i1}^{1}$ and
$H_{i2}^{1}$,and zero values for $\Pr\left[H_{i1}^{1}\neq H_{i2}^{1}\right],\Pr\left[H_{1j}^{2}\neq H_{2j}^{2}\right]$
are as acceptable a priori as any other. Analogous considerations
apply to $\left(H_{1j}^{2},H_{2j}^{2}\right)$ and $\left(B_{1j},B_{2j}\right)$. 

Our approach consists in replacing (\ref{eq:JDCassumption}) with
more general
\begin{equation}
\begin{array}{c}
\Pr\left[H_{i1}^{1}\neq H_{i2}^{1}\right]=2\varepsilon_{i}^{1}\in\left[0,1\right],\\
\Pr\left[H_{1j}^{2}\neq H_{2j}^{2}\right]=2\varepsilon_{j}^{2}\in\left[0,1\right],
\end{array}\label{eq:connections}
\end{equation}
and characterizing the dependence of $\left(A,B\right)$ on $\left(\alpha,\beta\right)$
by properties of the set of all 4-vectors $\varepsilon=\left(\varepsilon_{1}^{1},\varepsilon_{2}^{1},\varepsilon_{1}^{2},\varepsilon_{2}^{2}\right)$
that are compatible with or imply certain constraints imposed on the
vectors $p=\left(p_{11},p_{12},p_{21},p_{22}\right)$. Having adopted
a particular diagram of input-output correspondences (in our case,
$A\leftarrow\alpha,B\leftarrow\beta$), we can also say that these
sets of $\varepsilon$ characterize the contextual role of $\alpha,\beta$
for $B$ and $A$, respectively. 

We call $\varepsilon$ a vector of \emph{connection probabilities}.
The connection probabilities are of a principally non-empirical nature:
they are joint probabilities of events that can never co-occur. By
contrast, due to (\ref{eq:empirical joints}) the components of $p$
are joint probabilities of events that do co-occur, and by observing
these co-occurrences the probabilities in $p$ can be estimated. To
emphasize this distinction we refer to $p$ as a vector of \emph{empirical
probabilities}.

To distinguish our approach from other forms and meanings of probabilistic
contextualism, e.g., {[}\emph{\ref{enu:F.-Laudisa,-Contextualism},\ref{enu:A.-Yu.-Khrennikov,},\ref{enu:A.-Yu.-Khrennikov,-1}}{]},
we dub it the \emph{``all-possible-couplings''} approach. The term
``coupling'' refers to imposing a joint distribution (say, that
of $H_{11}^{1},H_{12}^{1}$) on random variables that otherwise are
not jointly distributed ($A_{11}$ and $A_{12}$). For a rigorous
and general discussion of couplings and connections see Appendix A.

\section{Extended Linear Feasibility Polytope (ELFP)}

ELFP is the set of all possible $\left(p,\varepsilon\right)$ for
which there exists a vector $H$ in (\ref{eq:extended H}) with jointly
distributed components $H_{ij}^{k}$ such that (\ref{eq:empirical joints})
holds, and, in accordance with (\ref{eq:connections}),
\begin{equation}
\Pr\left[H_{i1}^{1}=+1,H_{i2}^{1}=+1\right]=\varepsilon_{i}^{1},\;\Pr\left[H_{1j}^{2}=+1,H_{2j}^{2}=+1\right]=\varepsilon_{j}^{2},
\end{equation}
for all $i,j$. The existence of such an $H$ means the existence
of a probability vector $Q$ consisting of the $2^{8}$ joint probabilities
\begin{equation}
\Pr\left[H_{11}^{1}=h_{11}^{1},H_{11}^{2}=h_{11}^{2},\ldots,H_{22}^{2}=h_{22}^{2}\right],
\end{equation}
$h_{ij}^{1},h_{ij}^{2}\in\left\{ +1,-1\right\} $. Let $P$ denote
the $2{}^{5}\textnormal{-}$component vector consisting of 2$^{4}$
empirical probabilities 
\begin{equation}
\Pr\left[A_{ij}=a_{ij},B_{ij}=b_{ij}\right]
\end{equation}
and 2$^{4}$ connection probabilities 
\begin{equation}
\Pr\left[A_{i1}=a_{i1},A_{i2}=a_{i2}\right],\;\Pr\left[B_{1j}=b_{1j},B_{2j}=b_{2j}\right],
\end{equation}
$a_{ij},b_{ij}\in\left\{ +1,-1\right\} $. 

\textcolor{black}{Define a $2^{5}\times2^{8}$ Boolean matrix $M$
whose rows are enumerated in accordance with components of $P$ (i.e.,
by equalities $\left[A_{ij}=a_{ij},B_{ij}=b_{ij}\right]$, $\left[A_{i1}=a_{i1},A_{i2}=a_{i2}\right]$,
or $\left[B_{1j}=b_{1j},B_{2j}=b_{2j}\right]$) and columns in accordance
with components of $Q$ (i.e., by equalities $\left[H_{11}^{1}=h_{11}^{1},H_{11}^{2}=h_{11}^{2},\ldots,H_{22}^{2}=h_{22}^{2}\right]$).
An entry of $M$ contains 1 if and only if the corresponding random
variables in the enumerations of its row and its column have the same
values: e.g., if a row is enumerated by $\left[B_{12}=b_{12},B_{22}=b_{22}\right]$
and a column by $\left[H_{11}^{1}=h_{11}^{1},\ldots,H_{12}^{2}=h_{12}^{2},\ldots,H_{22}^{2}=h_{22}^{2}\right]$,
then their intersection contains 1 if and only if $h_{12}^{2}=b_{12},h_{22}^{2}=b_{22}$.}

It is easy to see that $H$ exists if and only if 
\begin{equation}
MQ=P
\end{equation}
for some vector $Q\ge0$ (componentwise) of probabilities. The vectors
$P$ for which such a $Q$ exists are exactly those within the polytope
whose vertices are the columns of the matrix $M$. The term ELFP is
due to this construction extending that of the linear feasibility
test in {[}\emph{\ref{enu:E.N.-Dzhafarov,-&2012}}{]}. This test,
among other applications, is the most general way of extending the
Bell/CHSH criterion to an arbitrary number of particles, spin axes,
and spin quantum numbers {[}\emph{\ref{enu:E.N.-Dzhafarov,-&2012},\ref{enu:E.N.-Dzhafarov,-in press 2},\ref{enu:I.-Pitowsky,-Quantum}-\ref{enu:R.M.-Basoalto,-&}}{]}.
Its application to binary inputs/outputs (not necessarily with equiprobable
outcomes) is shown in Appendix \ref{sub:Derivation-of-theFine}.

To describe ELFP by inequalities on $\left(p,\varepsilon\right)$,
we introduce the 16-component sets
\begin{equation}
\begin{array}{l}
\mathrm{\mathsf{S}}p=\left\{ \begin{array}{l}
\pm\left(p_{11}-\nicefrac{1}{4}\right)\pm\left(p_{12}-\nicefrac{1}{4}\right)\\
\pm\left(p_{21}-\nicefrac{1}{4}\right)\pm\left(p_{22}-\nicefrac{1}{4}\right):\\
\textnormal{each }\pm\textnormal{ is }+\textnormal{ or }-
\end{array}\right\} ,\\
\\
\mathrm{\mathsf{S}}\varepsilon=\left\{ \begin{array}{l}
\pm\left(\varepsilon_{1}^{1}-\nicefrac{1}{4}\right)\pm\left(\varepsilon_{1}^{2}-\nicefrac{1}{4}\right)\\
\pm\left(\varepsilon_{2}^{1}-\nicefrac{1}{4}\right)\pm\left(\varepsilon_{2}^{2}-\nicefrac{1}{4}\right):\\
\textnormal{each }\pm\textnormal{ is }+\textnormal{ or }-
\end{array}\right\} .
\end{array}
\end{equation}
$\mathrm{\mathsf{S}}_{0}p$ and $\mathrm{\mathsf{S}}_{1}p$ denote
the subsets of $\mathrm{\mathsf{S}}p$ with, respectively, even (0,2,
or 4) and odd (1 or 3) number of $+$ signs; $\mathrm{\mathsf{S}}_{0}\varepsilon$
and $\mathrm{\mathsf{S}}_{1}\varepsilon$ are defined analogously.
ELFP is described by
\begin{equation}
\max\left(\max\mathsf{S}_{0}p+\max\mathsf{S}_{1}\varepsilon,\;\max\mathsf{S}_{1}p+\max\mathsf{S}_{0}\varepsilon\right)\le\nicefrac{3}{2}\label{eq:ELFT}
\end{equation}
(see Appendix \ref{sub:Computations-for-ELFP}).

\section{All, Fit, Force, and Equi sets }

Let $\mathrm{constr}\left(p\right)$ denote any constraint (e.g.,
inequalities) imposed on $p$. Our approach consists in characterizing
this constraint by solving the following four problems:
\begin{enumerate}
\item Find the set $\mathrm{All{}_{constr}}$ of all $\left(p,\varepsilon\right)\in\left[0,\nicefrac{1}{2}\right]^{8}$
with $p$ subject to $\mathrm{constr}\left(p\right)$:
\begin{equation}
\left(p,\varepsilon\right)\in\mathrm{All_{constr}}\Longleftrightarrow\left(\mathrm{constr}(p)\textnormal{ and }(p,\varepsilon)\in\mathrm{ELFP}\right).
\end{equation}

\item Find the set $\mathrm{Fit_{constr}}$ of connection vectors $\varepsilon\in\left[0,\nicefrac{1}{2}\right]^{4}$
that fit (are compatible with) all empirical probability vectors $p$
satisfying $\mathrm{constr}$:
\begin{equation}
\varepsilon\in\mathrm{Fit}_{\mathrm{const}}\Longleftrightarrow(\mathrm{constr}(p)\Longrightarrow(p,\varepsilon)\in\mathrm{ELFP}).
\end{equation}

\item Find the set $\mathrm{Force}_{\mathrm{constr}}$ of $\varepsilon\in\left[0,\nicefrac{1}{2}\right]^{4}$
that force all compatible empirical probability vectors $p$ to satisfy
$\mathrm{constr}$:
\begin{equation}
\varepsilon\in\mathrm{Force}_{\mathrm{constr}}\Longleftrightarrow((p,\varepsilon)\in\mathrm{ELFP}\mathrm{\Longrightarrow constr}(p)).
\end{equation}

\item Find the set $\mathrm{Equi}_{\mathrm{constr}}$ of $\varepsilon\in\left[0,\nicefrac{1}{2}\right]^{4}$
for which an empirical probability vector $p$ satisfies $\mathrm{constr}$
if and only if $\left(p,\varepsilon\right)$ is in the ELFP set:
\begin{equation}
\varepsilon\in\mathrm{Equi}_{\mathrm{constr}}\Longleftrightarrow(\mathrm{constr}(p)\Longleftrightarrow(p,\varepsilon)\in\mathrm{ELFP}).
\end{equation}
Clearly, $\mathrm{Equi}_{\mathrm{constr}}=\mathrm{Force}_{\mathrm{constr}}\cap\mathrm{Fit}_{\mathrm{constr}}.$
\end{enumerate}
To illustrate, we focus on the following four benchmark constraints.
The no-constraint, or ``complete chaos'' situation is given by 
\begin{equation}
\mathrm{chaos(p)}\Longleftrightarrow p\in\left[0,\nicefrac{1}{2}\right]^{4},
\end{equation}
equivalent to (\ref{eq:chaos}) . The quantum mechanical constraint
is given by 
\begin{equation}
\mathrm{quant}(p)\Longleftrightarrow\max\mathsf{S}_{1}p\leq\nicefrac{\sqrt{2}}{2},
\end{equation}
equivalent to (\ref{eq:Cirelson}) . The ``classical'' constraint
is given by 
\begin{equation}
\mathrm{class}(p)\Longleftrightarrow\max\mathsf{S}_{1}p\leq\nicefrac{1}{2},
\end{equation}
equivalent to the Bell/CHSH inequalities (\ref{eq:Fine}). Finally,
we consider the constraint 
\begin{equation}
\mathrm{fix}\left(p\right)\Longleftrightarrow p=\textnormal{specific vector}.
\end{equation}
For all constraints except for $\mathrm{fix}\left(p\right)$ the sets
All, Fit, Force, and Equi are as shown in Table 1 (for derivations
see Appendix \ref{sub:Computations-for-3,}).

\begin{center}
\begin{table}[h]
\caption{Characterizations of the sets of four different types (columns) subject
to three constrains (rows).\label{tab:Characterization} }

\centering{}%
\begin{tabular}{|l|c|c|c|c|}
\hline 
 & All $\left(p,\varepsilon\right)$ & Fit $\left(\varepsilon\right)$ & Force $\left(\varepsilon\right)$ & Equi $\left(\varepsilon\right)$\tabularnewline
\hline 
\hline 
chaos & $\in$ ELFP & $\max\mathsf{S}\varepsilon\le\nicefrac{1}{2}$ & $\varepsilon\in\left[0,\frac{1}{2}\right]^{4}$ & $\max\mathsf{S}\varepsilon\le\nicefrac{1}{2}$\tabularnewline
\hline 
quant & $\begin{array}{c}
\max\mathsf{S}_{1}p\leq\nicefrac{\sqrt{2}}{2}\\
\&\\
\left(p,\varepsilon\right)\in\textnormal{ ELFP}
\end{array}$ & $\begin{array}{l}
\max\mathsf{S}_{0}\varepsilon\le\frac{3-\sqrt{2}}{2},\\
\max\mathsf{S}_{1}\varepsilon\le\nicefrac{1}{2}
\end{array}$ & $\max\mathsf{S}_{0}\varepsilon\geq\frac{3-\sqrt{2}}{2}$ & $\begin{array}{l}
\frac{3-\sqrt{2}}{2}\in\mathsf{S}_{0}\varepsilon,\\
\max\mathsf{S}_{1}\varepsilon\leq\nicefrac{1}{2}
\end{array}$\tabularnewline
\hline 
class & $\begin{array}{c}
\max\mathsf{S}_{1}p\leq\nicefrac{1}{2}\\
\&\\
\left(p,\varepsilon\right)\in\textnormal{ ELFP}
\end{array}$ & $\max\mathsf{S}_{1}\varepsilon\le\nicefrac{1}{2}$ & $1\in\mathsf{S}_{0}\varepsilon$ & $1\in\mathsf{S}_{0}\varepsilon$\tabularnewline
\hline 
\end{tabular}
\end{table}

\par\end{center}

Thus, $\mathrm{Fit_{chaos}}$ is the set of all $\varepsilon$ such
that $\max\mathsf{S}\varepsilon\le\nicefrac{1}{2}$: if an $\varepsilon$
is in this set, then any $p$ (with no constraints) is compatible
with it. $\mathrm{Force}_{\mathrm{quant}}$ is characterized by $\max\mathsf{S}_{0}\varepsilon\geq\frac{3-\sqrt{2}}{2}$:
if an $\varepsilon$ is in this set, then all compatible with it $p$
satisfy $\mathrm{quant}(p)$. $\mathrm{Equi}_{\mathrm{class}}$ is
the set of all $\varepsilon$ such that $\mathsf{S}_{0}\varepsilon$
contains $1$: for any such an $\varepsilon$, a $p$ is compatible
with it if and only if it satisfies $\mathrm{class}\left(p\right)$.

\begin{figure}
\begin{centering}
\includegraphics[bb=20bp 10bp 760bp 582bp,clip,scale=0.4]{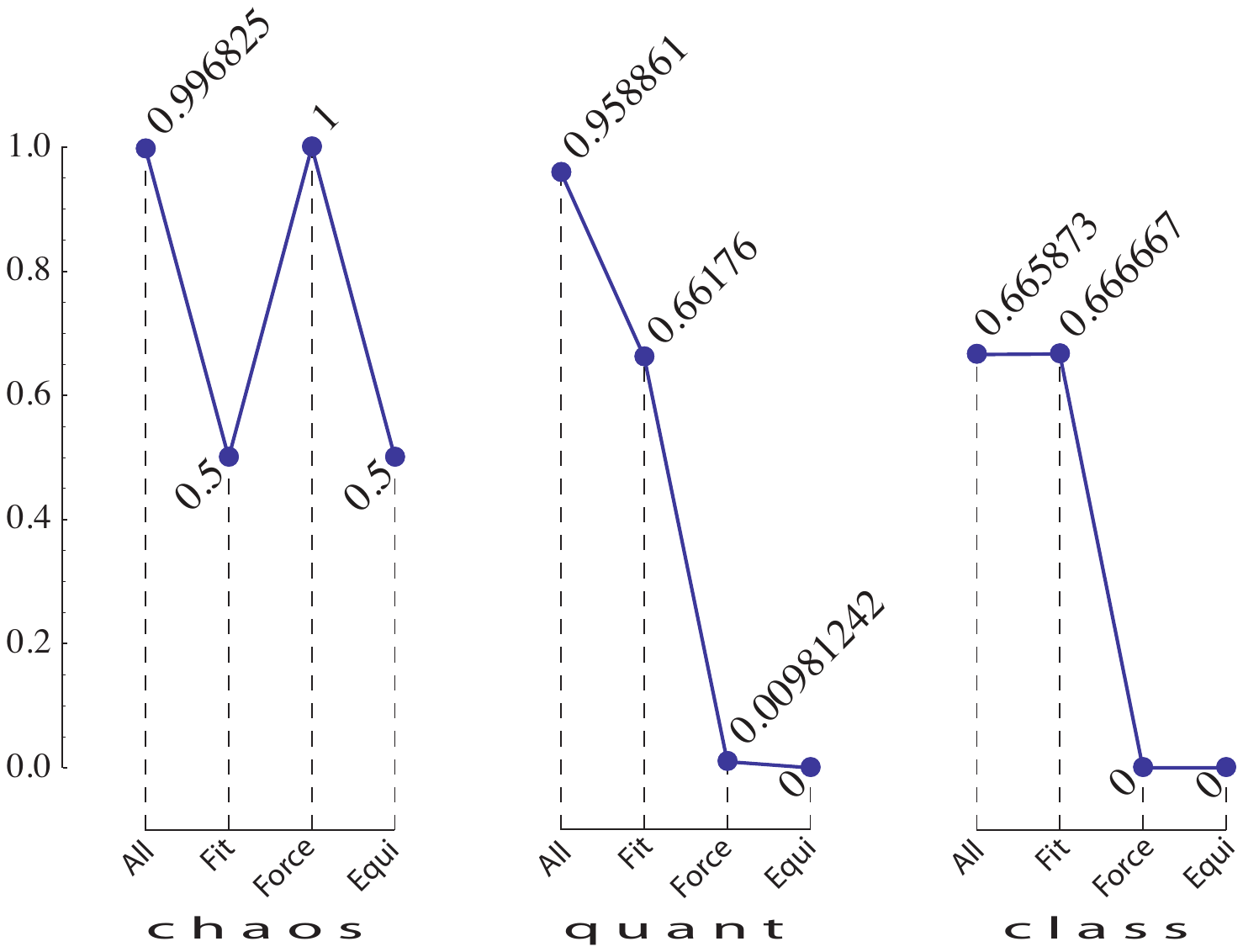}\includegraphics[bb=0bp 100bp 584bp 344bp,scale=0.4]{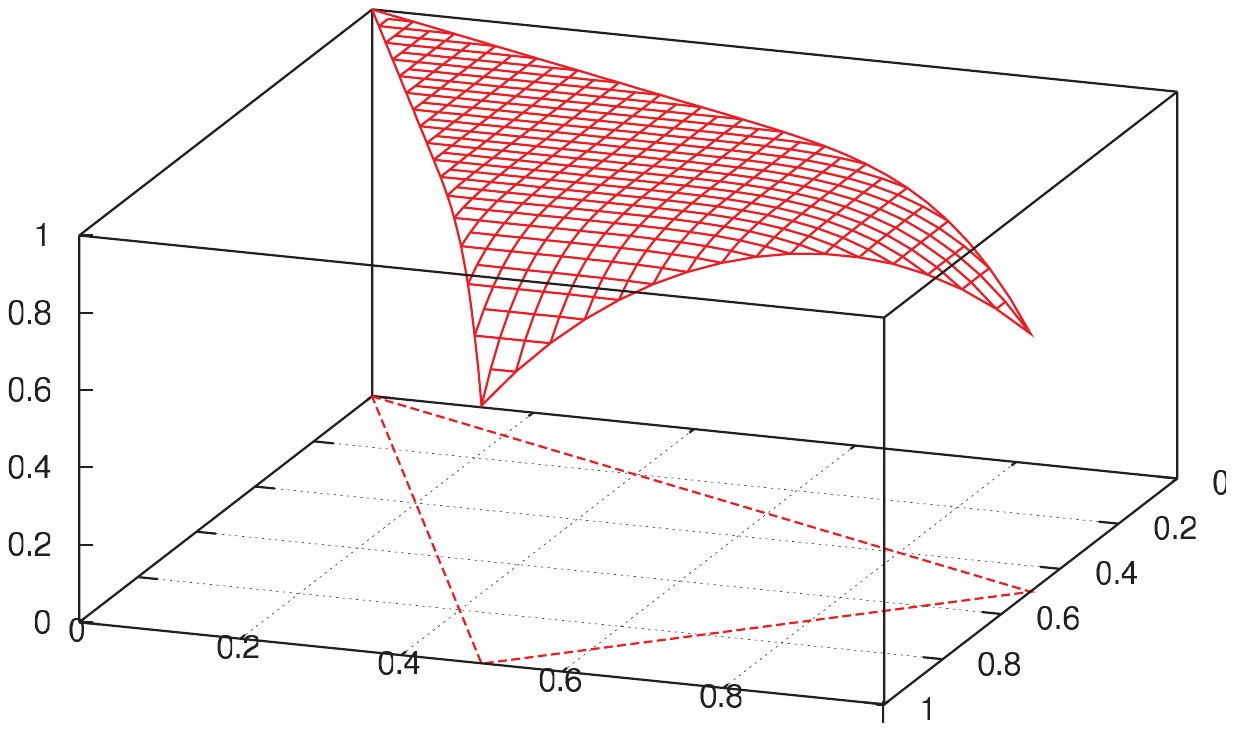}
\par\end{centering}

\caption{Left: Profiles $\text{Vol}^{8}\left(\mathrm{All}_{\mathrm{constr}}\right)\rightarrow\text{Vol}^{4}\left(\mathrm{Fit}_{\mathrm{constr}}\right)\rightarrow\text{Vol}^{4}\left(\mathrm{Force}_{\mathrm{constr}}\right)\rightarrow\text{Vol}^{4}\left(\mathrm{Equi}_{\mathrm{constr}}\right)$
for constraints chaos, quant, and class. Right: $\mathrm{Vol}^{4}\left(\mathrm{Fit}_{\mathrm{fix}\left(p\right)}\right)$
as a function of $x=\max\mathsf{S}_{0}p$ and $y=\max\mathsf{S}_{1}p$.
The possible $\left(x,y\right)\textnormal{-}$pairs form the triangle
$\left((0,0),(\nicefrac{1}{2},1),(1,\nicefrac{1}{2})\right)$, and
$\mathrm{Vol}^{4}\left(\mathrm{Fit}_{\mathrm{fix}\left(p\right)}\right)=1+\frac{\rho\left(x\right)}{3}\left(-1+8x-24x^{2}+32x^{3}-16x^{4}\right)+\frac{\rho\left(y\right)}{3}\left(-1+8y-24y^{2}+32y^{3}-16y^{4}\right)$,
where $\rho\left(z\right)=1$ if $z\geq\nicefrac{1}{2}$ and $\rho\left(z\right)=0$
otherwise. }

\end{figure}

For each of these sets we compute $\mathrm{Vol}^{d}$, its volume
normalized by that of $\left[0,\nicefrac{1}{2}\right]^{d}$, with
$d$ being the dimensionality of the set (Fig.$\:$2, left). Thus,
the defining property of $\mathrm{Force}_{\mathrm{class}}$, $1\in\mathsf{S}_{0}\varepsilon$,
is satisfied if and only if either all $\varepsilon_{i}^{k}$ are
$0$, or they all are $\nicefrac{1}{2}$, or two of them are $0$
and two $\nicefrac{1}{2}$. Hence $\text{Vol}^{4}\left(\mathrm{Force}_{\mathrm{class}}\right)=0$.
For nonzero volumes, the derivation is described in Appendix \ref{sub:Computations-for-3,}.
Each panel of Fig.$\:$2, left, can be viewed as a ``profile'' of
the corresponding constraint. Each of the first three volumes in a
panel can be viewed as characterizing the ``strictness'' of a constraint,
in three different meanings. The intuition of a stricter constraint
is that it corresponds to a smaller $\text{Vol}^{8}\left(\mathrm{All}_{\mathrm{constr}}\right)$,
larger $\text{Vol}^{4}\left(\mathrm{Fit}_{\mathrm{constr}}\right)$,
and smaller $\text{Vol}^{4}\left(\mathrm{Force}_{\mathrm{constr}}\right)$.
Characterizing constraints imposed on empirical probabilities by multidimensional
volumes is not a new idea {[}\ref{enu:Cabello-A-(2005)}{]}, but our
computations are different: they are aimed at sets of nonempirical
connection probabilities in relation to constraints imposed on empirical
probabilities.

The constraint $\mathrm{fix}\left(p\right)$ has to be handled separately.
Clearly, $\text{Vol}^{8}\left(\mathrm{All}_{\mathrm{\mathrm{fix}\left(p\right)}}\right)=0$.
$\mathrm{Fit}_{\mathrm{fix}\left(p\right)}$ is described by

\begin{equation}
\begin{array}{l}
\max\mathsf{S}_{1}\varepsilon\le\nicefrac{3}{2}-\max\mathsf{S}_{0}p,\\
\max\mathsf{S}_{0}\varepsilon\le\nicefrac{3}{2}-\max\mathsf{S}_{1}p,
\end{array}\label{eq:fixed}
\end{equation}
and $\mathrm{Vol}^{4}\left(\mathrm{Fit}_{\mathrm{fix}\left(p\right)}\right)$
is a polynomial function of $\max\mathsf{S}_{0}p$ and $\max\mathsf{S}_{1}p$,
these two quantities forming the triangle $\left((0,0),(\nicefrac{1}{2},1),(1,\nicefrac{1}{2})\right)$.
The polynomial and its values are shown in Fig.$\:$2, right (see
Appendix \ref{sub:Computations-for-1}, for computational details).
$\mathrm{Force}_{\mathrm{fix}\left(p\right)}$ is clearly empty, hence
so is $\mathrm{Equi}_{\mathrm{fix}\left(p\right)}$.

\section{Conclusion}

The essence of the proposed mathematical framework is as follows.
We consider all possible couplings for empirically observed vectors
of random outputs. In the case of two binary inputs/outputs these
vectors are pairs 
\[
\left(A_{11},B_{11}\right),\left(A_{12},B_{12}\right),\left(A_{21},B_{21}\right),\left(A_{22},B_{22}\right),
\]
the couplings $H$ for them have the form (\ref{eq:extended H}),
with the coupling relation (\ref{eq:empirical joints}). We assume
that the joint distributions (in our case described by pairwise joint
probabilities) of the empirically observed $\left(A_{ij},B_{ij}\right)$
are subject to a certain constraint, given to us by substantive considerations
outside the scope of our approach: for instance, if a system consists
of entangled particles, a constraint, say (\ref{eq:Cirelson}), is
derived from the quantum theory. Due to (\ref{eq:empirical joints}),
the constraint is imposed on
\begin{equation}
\left(H_{11}^{1},H_{11}^{2}\right),\left(H_{12}^{1},H_{12}^{2}\right),\left(H_{21}^{1},H_{21}^{2}\right),\left(H_{22}^{1},H_{22}^{2}\right).\label{eq:observable part}
\end{equation}
We investigate then the unobservable ``connections'', the subvectors
of the components of $H$ that correspond to outputs obtained at mutually
exclusive values of the inputs (i.e., never co-occurring). In our
case these are the pairs 
\begin{equation}
\left(H_{11}^{1},H_{12}^{1}\right),\left(H_{21}^{1},H_{22}^{1}\right),\left(H_{11}^{2},H_{21}^{2}\right),\left(H_{12}^{2},H_{22}^{2}\right)\label{eq:unobservable part}
\end{equation}
corresponding to, respectively,
\begin{equation}
\left(A_{11},A_{12}\right),\left(A_{21},A_{22}\right),\left(B_{11},B_{21}\right),\left(B_{12},B_{22}\right).
\end{equation}
 We then characterize the constraint imposed on the empirical pairs
(\ref{eq:observable part}) by describing the ``fitting'' or ``forcing''
(or both ``fitting and forcing'') distributions of the unobservable
connections (\ref{eq:unobservable part}). By fitting distributions
of (\ref{eq:unobservable part}) we mean those that are compatible
with any (\ref{eq:observable part}) subject to the constraint in
question, the compatibility meaning that all these eight pairs can
be embedded into a single $H$ (with jointly distributed components).
By forcing distributions of (\ref{eq:unobservable part}) we mean
those that are compatible with (\ref{eq:observable part}) only if
the latter are subject to the given constraint. 

The value of this approach is in providing a unified language for
speaking of probabilistic contextuality. At the cost of greater computational
complexity but with no conceptual complications the computations involved
in our demonstration of the all-possible-couplings approach can be
extended to more general cases: arbitrary marginal probabilities (satisfying
marginal selectivity), nonlinear constraints, and greater numbers
of inputs, outputs, and their possible values. The language for a
completely general theory, involving unrestricted (not necessarily
finite) sets of inputs, outputs, and their values, is presented in
Appendix A.

\FloatBarrier

\pagebreak{}

\renewcommand{\thesection}{A}

\setcounter{equation}{0}

\setcounter{thm}{0}

\section*{\label{sub:Coupling-approach}APPENDIX A: All-possible-couplings
approach on the general level}

We show here how the approach presented in the main text generalizes
to arbitrary sets of inputs and random outputs. We use the term \emph{sequence}
to refer to any indexed family (a function from an index set into
a set), with index sets not necessarily countable.\emph{ }We present
sequences in the form $\left(x^{y}:y\in Y\right)$, $\left(x_{z}:z\in Z\right)$,
or $\left(x_{z}^{y}:y\in Y,z\in Z\right)$. A random variable is understood
most broadly, as a measurable mapping between any two probability
spaces. In particular, any sequence of jointly distributed random
variables is a random variable. For brevity, we omit an explicit presentation
of probability spaces and distributions. In all other respects the
notation and terminology closely follow {[}\emph{\ref{enu:E.N.-Dzhafarov,-&in press},\ref{enu:E.N.-Dzhafarov,-in press 2}}{]}. 

An \emph{input} is a set of elements called \emph{input values}. Let
$\alpha=\left(\alpha^{k}:k\in K\right)$ be a sequence of inputs.
A \emph{treatment} is a sequence $\phi=(x^{k}:k\in K)$ that belongs
to a nonempty set $\Phi\subset\prod_{k\in K}\alpha^{k}$ (so that
$x^{k}\in\alpha^{k}$ for all $k\in K$). If $\phi\in\Phi$, $k\in K$,
and $I\subset K$, then $\phi\left(k\right)=x^{k}\in\alpha^{k}$ and
$\phi|I$ is the restriction of $\phi$ to $I$, i.e., the sequence
$(x^{k}:k\in I)$.

An \emph{output} is a random variable. Let $\left(A_{\phi}^{k}:k\in K,\phi\in\Phi\right)$
be a sequence of outputs such that 
\begin{enumerate}
\item $A_{\phi}=\left(A_{\phi}^{k}:k\in K\right)$ is a random variable
for every $\phi\in\Phi$, i.e., the random variables $A_{\phi}^{k}$
across all possible $k$ possess a joint distribution;
\item if $\phi,\phi'\in\Phi$, $I\subset K$, and $\phi|I=\phi'|I$, then
$\left(A_{\phi}^{k}:k\in I\right)\sim\left(A_{\phi'}^{k}:k\in I\right)$.
\end{enumerate}
Property 2 is \emph{(complete) marginal selectivity} {[}\emph{\ref{enu:Dzhafarov,-E.N.-(2003).}}{]}.
$A_{\phi}$ is called an \emph{empirical random variable}, and $A=\left(A_{\phi}:\phi\in\Phi\right)$
is the \emph{sequence of empirical random variables}.
\begin{rem}
The interpretation is that for every $\phi$, each $\alpha^{k}$ may
``directly'' influence $A_{\phi}^{k}$ but no other output in $A_{\phi}$.
The fact that inputs in $\alpha=\left(\alpha^{k}:k\in K\right)$ and
outputs in an \emph{empirical random variable} $A_{\phi}=\left(A_{\phi}^{k}:k\in K\right)$
are in a bijective correspondence is not restrictive: this can always
be achieved by an appropriate grouping of inputs and (re)definition
of treatments $\phi$ {[}\emph{\ref{enu:E.N.-Dzhafarov,-&2012}}{]}. 
\end{rem}

\begin{rem}
The case considered in the main text corresponds to $K=\left\{ 1,2\right\} $,
\begin{equation}
\alpha=\left(\alpha^{1},\alpha^{2}\right)\textnormal{ with }\alpha^{k}=\left\{ \alpha_{1}^{k},\alpha_{2}^{k}\right\} \textnormal{ for }k\in\left\{ 1,2\right\} ,
\end{equation}
\begin{equation}
\begin{array}{l}
\Phi=\left\{ \phi_{11},\phi_{12},\phi_{21},\phi_{22}\right\} \\
\textnormal{ with }\phi_{ij}=\left(\alpha_{i}^{1},\alpha_{j}^{2}\right)\textnormal{ for }i,j\in\left\{ 1,2\right\} ,
\end{array}
\end{equation}
and (abbreviating $A_{\phi_{ij}}$ as $A_{ij}$ and $A_{\phi_{ij}}^{k}$
as $A_{ij}^{k}$)
\begin{equation}
\begin{array}{l}
A=\left(A_{11},A_{12},A_{21},A_{22}\right),\\
\textnormal{ with }A_{ij}=\left(A_{i}^{1},A_{j}^{2}\right)\textnormal{ for }i,j\in\left\{ 1,2\right\} ,
\end{array}
\end{equation}
 where each $A_{ij}^{k}$ is a binary random variable with $\Pr\left[A_{ij}^{k}=a_{1}^{k}\right]=\Pr\left[A_{ij}^{k}=a_{2}^{k}\right]=\nicefrac{1}{2}$. 

\medskip{}

\end{rem}
Given a sequence of \emph{empirical random variables} $A=\left(A_{\phi}:\phi\in\Phi\right)$,
a sequence of random variables 
\begin{equation}
C_{A}=\left(C_{\tau}^{I}:\tau\in\prod_{k\in I}\alpha^{k},I\in2^{K}-\left\{ \emptyset,K\right\} \right)
\end{equation}
(not necessarily jointly distributed) is called a \emph{connecting
set} for $A$ if each $C_{\tau}^{I}$ is a coupling for 
\begin{equation}
A_{\tau}^{I}=\left(A_{\phi}^{I}:\phi\in\Phi,\phi|I=\tau\right),
\end{equation}
where $A_{\phi}^{I}=\left(A_{\phi}^{k}:k\in I\right)$. This means
that $C_{\tau}^{I}$ is a random variable of the form 
\begin{equation}
C_{\tau}^{I}=\left(C_{\tau,\phi}^{I}:\phi\in\Phi,\phi|I=\tau\right)
\end{equation}
with 
\begin{equation}
C_{\tau,\phi}^{I}\sim A_{\phi}^{I}
\end{equation}
for all $\phi\in\Phi$ such that $\phi|I=\tau$. $C_{\tau}^{I}$ is
called an \emph{$\left(I,\tau\right)\textnormal{-}$connection}. \textcolor{black}{The
indexation in $C_{\tau,\phi}^{I}$ is to ensure that if $\left(I,\tau\right)\not=\left(I',\tau'\right)$,
then $C_{\tau}^{I}$ and $C_{\tau'}^{I'}$ are stochastically unrelated.}
An \emph{identity $\left(I,\tau\right)\textnormal{-}$connection}
$C_{\tau}^{I}$ is one with $\Pr\left[C_{\tau,\phi}^{I}=C_{\tau,\phi'}^{I}\right]=1$
for any $\phi,\phi'\in\Phi$. Formally, $A$ itself can be viewed
as a connection $C_{\emptyset}^{\emptyset}$, but it is preferable
to keep $A$ separate by not allowing $I=\emptyset$, as $A$ is the
only empirically observable part of the construction. 
\begin{rem}
It is generally convenient not to distinguish identically distributed
connections. By abuse of language, the distribution of $C_{\tau}^{I}$
(or some characterization thereof) can also be called \emph{$\left(I,\tau\right)\textnormal{-}$}connection.
We use this language in the main text when we represent\textcolor{red}{{}
}\textcolor{black}{\emph{$\left(\left\{ k\right\} ,k\mapsto\alpha_{i}^{k}\right)\textnormal{-}$}}\textcolor{black}{connections}
(without introducing them explicitly) by probabilities $\varepsilon_{i}^{k}$
and call $\varepsilon$ a connection vector. See Remark \ref{rem:connections}.

\medskip{}

\end{rem}
A jointly distributed sequence
\begin{equation}
H=\left(H_{\phi}^{k}:k\in K,\phi\in\Phi\right)
\end{equation}
is called an \emph{Extended Joint Distribution Sequence} (EJDS) for
$\left(A,C_{A}\right)$ if for any $I\in2^{K}-\left\{ \emptyset,K\right\} $
and any $\tau\in\prod_{k\in I}\alpha^{k}$,
\begin{equation}
H_{\tau}^{I}=\left(H_{\phi}^{I}:\phi\in\Phi,\phi|I=\tau\right)\sim C_{\tau}^{I},
\end{equation}
where $H_{\phi}^{I}=\left(H_{\phi}^{k}:k\in I\right)$, and 
\begin{equation}
H_{\phi}^{K}=\left(H_{\phi}^{k}:k\in K\right)\sim A_{\phi}
\end{equation}
for any $\phi\in\Phi$.
\begin{rem}
\label{rem:connections}For the case considered in the main text,
a connecting set for $A$ is (conveniently replacing $C_{\phi_{ij}}^{\left\{ k\right\} }$
, $C_{\phi_{ij}|\left\{ 1\right\} }^{\left\{ 1\right\} }$, and $C_{\phi_{ij}|\left\{ 2\right\} }^{\left\{ 2\right\} }$
with $C_{ij}^{k}$ , $C_{i}^{1}$, and $C_{j}^{2}$, respectively)
\begin{equation}
\begin{array}{l}
C_{A}=\left(C_{1}^{1},C_{2}^{1},C_{1}^{2},C_{2}^{2}\right)\\
\textnormal{ with }C_{i}^{1}=\left(C_{i,i1}^{1},C_{i,i2}^{1}\right)\textnormal{ and }C_{j}^{2}=\left(C_{j,1j}^{2},C_{j,2j}^{2}\right),
\end{array}
\end{equation}
such that
\begin{equation}
C_{i,ij}^{1}\sim A_{ij}^{1},\; C_{j,ij}^{2}\sim A_{ij}^{2}\textnormal{ for }i,j\in\left\{ 1,2\right\} .
\end{equation}
An EJDS for $\left(A,C_{A}\right)$ is a random variable (using analogous
abbreviations)
\begin{equation}
H=\left(H_{11}^{1},H_{11}^{2},H_{12}^{1},H_{12}^{2},H_{21}^{1},H_{21}^{2},H_{22}^{1},H_{22}^{2}\right)
\end{equation}
such that
\begin{equation}
\left(H_{i1}^{1},H_{i2}^{1}\right)\sim C_{i}^{1},\;\left(H_{1j}^{2},H_{2j}^{2}\right)\sim C_{j}^{2}
\end{equation}
and
\begin{equation}
H_{ij}^{12}=\left(H_{ij}^{1},H_{ij}^{2}\right)\sim A_{ij}=\left(A_{ij}^{1},A_{ij}^{2}\right)\textnormal{ for }i,j\in\left\{ 1,2\right\} .
\end{equation}
In the main text each $C_{i}^{k}$ is represented by $\varepsilon_{i}^{k}$
and each $H_{ij}^{12}$ by $p_{ij}$. 

\medskip{}

\end{rem}
An EJDS for $\left(A,C_{A}\right)$ reduces to the Joint Distribution
Criterion set (JDC set) of the theory of selective influences {[}\emph{\ref{enu:Dzhafarov,-E.N.,-&2010}-\ref{enu:E.N.-Dzhafarov,-in press 2}}{]}
if all connections in $C_{A}$ are identity ones. Note that no connection
has an empirical meaning: for distinct $\phi,\phi'\in\Phi$, the variables
$A_{\phi}^{I}$ and $A_{\phi'}^{I}$ corresponding to $C_{\tau,\phi}^{I}$
and $C_{\tau,\phi'}^{I}$ do not have an empirically observable (or
theoretically privileged) pairing scheme. 

Let $X$ be any set whose elements are sequences of \emph{empirical
random variables} $A=\left(A_{\phi}:\phi\in\Phi\right)$. $X$ can
be viewed as the set of all possible \emph{empirical random variables}
satisfying certain constraints. We define the sets All$_{X}$, Fit$_{X}$,
Force$_{X}$, and Equi$_{X}$ as follows:
\begin{enumerate}
\item All$_{X}$ is the set of all pairs $\left(A,C_{A}\right)$ such that
\begin{equation}
A\in X\textnormal{ and there exists an EJDS }H\textnormal{ for }\left(A,C_{A}\right).
\end{equation}

\item Fit$_{X}$ is the set of all $C_{A}$ such that 
\begin{equation}
A\in X\Longrightarrow\textnormal{there exists an EJDS }H\textnormal{ for }\left(A,C_{A}\right).
\end{equation}

\item Force$_{X}$ is the set of all $C_{A}$ such that
\begin{equation}
\textnormal{there exists an EJDS }H\textnormal{ for }\left(A,C_{A}\right)\Longrightarrow A\in X.
\end{equation}

\item $\mathrm{Equi}_{X}=\mathrm{Force}_{X}\cap\mathrm{Fit}_{X},$ that
is, $C_{A}\in\mathrm{Equi}_{X}$ if and only if
\begin{equation}
A\in X\Longleftrightarrow\textnormal{there exists an EJDS }H\textnormal{ for }\left(A,C_{A}\right).
\end{equation}

\end{enumerate}
The all-possible-couplings approach in the general case consists in
characterizing any $X$ (interpreted as a type of contextuality or
determinism) by All$_{X}$, Fit$_{X}$, Force$_{X}$, and Equi$_{X}$.
A straightforward generalization of this approach that might be useful
in some applications is to replace $C_{A}$ in all definitions with
a subset of $C_{A}$, or several subsets of $C_{A}$ tried in turn.
Thus one might consider connections involving only particular $I\subset K$
(e.g., only singletons), or one might require that some of the connections
are identity ones. 

\renewcommand{\thesection}{B}

\renewcommand{\thesubsection}{B\arabic{subsection}}

\setcounter{subsection}{0}

\setcounter{equation}{0}

\setcounter{thm}{0}

\pagebreak{}

\section*{\label{sec:APPENDIX-B:-Technical}APPENDIX B: Technical Details}

\subsection{\label{sub:Derivation-of-theFine}Derivation of the Bell/CHSH bounds}

A representation (\ref{eq:JDC1})-(\ref{eq:JDC2}) exists if and only
if the 2$^{4}$ possible values $\left(h_{1}^{1},h_{2}^{1},h_{1}^{2},h_{2}^{2}\right)$
of $H$ ($h_{i}^{k}\in\left\{ +1,-1\right\} ,i,k\in\left\{ 1,2\right\} $)
can be assigned probabilities
\begin{equation}
p\left(h_{1}^{1},h_{2}^{1},h_{1}^{2},h_{2}^{2}\right)=\Pr\left[H_{1}^{1}=h_{1}^{1},H_{2}^{1}=h_{2}^{1},H_{1}^{2}=h_{1}^{2},H_{2}^{2}=h_{2}^{2}\right],
\end{equation}
so that, for all $a_{ij},b_{ij}\in\left\{ +1,-1\right\} $, $i,j\in\left\{ 1,2\right\} $,
\begin{equation}
\sum_{h_{1}^{1},h_{2}^{1},h_{1}^{2},h_{2}^{2}}\chi\left(h_{i}^{1}=a_{ij}\,\wedge\, h_{j}^{2}=b_{ij}\right)p\left(h_{1}^{1},h_{2}^{1},h_{1}^{2},h_{2}^{2}\right)=\Pr\left[A_{ij}=a_{ij},B_{ij}=b_{ij}\right],
\end{equation}
where $\chi\left(...\right)$ indicates the truth value (1 or 0) of
the statement within the parentheses. It is easy to see that this
system of linear equations can be written as
\begin{equation}
MQ=P,\label{eq:LFT}
\end{equation}
where $P$ is the 16-vector of probabilities $\Pr\left[A_{ij}=a_{ij},B_{ij}=b_{ij}\right]$
indexed (together with the columns of matrix $M$) by $\left(i,j,a_{ij},b_{ij}\right)$-values,
say, lexicographically; $Q$ is the 16-vector of unknown probabilities
$p\left(h_{1}^{1},h_{2}^{1},h_{1}^{2},h_{2}^{2}\right)$ indexed (together
with the rows of $M$) by $\left(h_{1}^{1},h_{2}^{1},h_{1}^{2},h_{2}^{2}\right)$-values
in some order; and the cells of $M$ indexed by $\left(\left(i,j,a_{ij},b_{ij}\right),\left(h_{1}^{1},h_{2}^{1},h_{1}^{2},h_{2}^{2}\right)\right)$
contain $\chi\left(h_{i}^{1}=a_{ij}\,\wedge\, h_{j}^{2}=b_{ij}\right)$.
We conclude that a representation (\ref{eq:JDC1})-(\ref{eq:JDC2})
exists if and only if 
\begin{equation}
B\left(M,P\right)=1,\label{eq:B(M,P)}
\end{equation}
where $B$$\left(M,P\right)$ is a Boolean function equal to 1 if
(\ref{eq:LFT}) has at least one solution with nonnegative components
of $Q$. It is easy to show (see {[}\ref{enu:E.N.-Dzhafarov,-&2012}{]}
for details) that solutions $Q$ of (\ref{eq:LFT}) always have the
property
\begin{equation}
\sum_{h_{1}^{1},h_{2}^{1},h_{1}^{2},h_{2}^{2}}p\left(h_{1}^{1},h_{2}^{1},h_{1}^{2},h_{2}^{2}\right)=1.
\end{equation}
It is known from the linear programming theory that $B\left(M,P\right)$
is always computable. A standard facet enumeration algorithm allows
one to obtain the system of all linear inequalities and equations
imposed on $P$ that are equivalent to (\ref{eq:B(M,P)}). This system
turns out to consist of the equalities (\ref{eq:marginal selectivity})
representing marginal selectivity, and inequalities that can be written
as
\begin{equation}
-2\leq\textnormal{E}_{ij}+\textnormal{E}_{i'j}+\textnormal{E}_{i'j'}-\textnormal{E}_{ij'}\leq2,
\end{equation}
where, in reference to (\ref{eq:matrix}), $\textnormal{E}_{ij}=p_{ij}+s_{ij}-q_{ij}-r_{ij}$
is the expected value of $A_{ij}B_{ij}$. When marginal probabilities
are all $\nicefrac{1}{2}$, these inequalities reduce to (\ref{eq:Fine}),
using $p_{ij}=\left(\textnormal{E}_{ij}+1\right)/4$.
\begin{rem}
It would be a mistake to consider this proof ``computer-assisted''
because it mentions a facet enumeration algorithm. The latter is merely
a long chain of trivial algebraic transformations, that can always
be written out \emph{in extenso} if needed.
\end{rem}

\subsection{\label{sub:Derivation-of-theCirelson}Derivation of the Cirel'son
bounds.}

The following is a modification of the derivation given in {[}\ref{enu:Landau-LJ-(1987)}{]}.
Let $a,a',b,b'$ be the Hermitian operators in complex Hilbert space
corresponding to, respectively, outputs $A_{1j},A_{2j},B_{i1},B_{i2}$
(where $i$ and $j$ are irrelevant, i.e., $a$ represents both $A_{11}$
and $A_{12}$, $b$ both $B_{11}$ and $B_{21}$, etc.). Denoting
by $\mathrm{E}$ expected value and by $\mathrm{Tr}$ trace, we have,
for any state (density operator) $W$,
\begin{equation}
\begin{array}{c}
4p_{11}-1=\mathrm{E}\left[A_{11}B_{11}\right]=\mathrm{Tr}\left(Wab\right),\\
4p_{12}-1=\mathrm{E}\left[A_{12}B_{12}\right]=\mathrm{Tr}\left(Wab'\right),\\
etc.
\end{array}\label{eq:p-Exp}
\end{equation}
where either of $a$ and $a'$ commutes with either of $b$ and $b'$.
Inequalities (\ref{eq:Cirelson}) to be demonstrated are equivalent
to
\begin{equation}
\begin{array}{c}
R_{1}=\left|\mathrm{Tr}\left(Wab\right)+\mathrm{Tr}\left(Wab'\right)+\mathrm{Tr}\left(Wa'b\right)-\mathrm{Tr}\left(Wa'b'\right)\right|=\left|\mathrm{Tr}\left(Ws_{1}\right)\right|\leq2\sqrt{2},\\
\begin{array}{c}
R_{2}=\left|\mathrm{Tr}\left(Wab\right)+\mathrm{Tr}\left(Wab'\right)-\mathrm{Tr}\left(Wa'b\right)+\mathrm{Tr}\left(Wa'b'\right)\right|=\left|\mathrm{Tr}\left(Ws_{2}\right)\right|\leq2\sqrt{2},\\
etc.
\end{array}
\end{array}\label{eq:Landau}
\end{equation}
where
\begin{equation}
\begin{array}{c}
s_{1}=ab+ab'+a'b-a'b'=a\left(b+b'\right)+a'\left(b-b'\right),\\
s_{2}=ab+ab'-a'b+a'b'=a\left(b+b'\right)-a'\left(b-b'\right),\\
etc.
\end{array}
\end{equation}
Since the values of the outputs, +1/-1, are the eigenvalues of the
corresponding operators, it can easily be seen (e.g., by spectral
decomposition, squaring, and then multiplication by an arbitrary vector)
that 
\begin{equation}
a^{2}=b^{2}=a'^{2}=b'^{2}=I,
\end{equation}
where $I$ is the identity operator. Using this we show by straightforward
if somewhat tedious algebra that
\begin{equation}
\begin{array}{c}
s_{1}^{2}=s_{4}^{2}=4I-\left(aa'-a'a\right)\left(bb'-b'b\right),\\
s_{2}^{2}=s_{3}^{2}=4I+\left(aa'-a'a\right)\left(bb'-b'b\right),
\end{array}
\end{equation}
whence, using the conventional notation for commutators, $\left[x,y\right]=xy-yx$,
\begin{equation}
\begin{array}{c}
\mathrm{Tr}\left(Ws_{1}^{2}\right)=\mathrm{Tr}\left(Ws_{4}^{2}\right)=4-\mathrm{Tr}\left(W\left[a,a'\right]\left[b,b'\right]\right),\\
\mathrm{Tr}\left(Ws_{2}^{2}\right)=\mathrm{Tr}\left(Ws_{3}^{2}\right)=4+\mathrm{Tr}\left(W\left[a,a'\right]\left[b,b'\right]\right).
\end{array}\label{eq:traces}
\end{equation}
For $k=1,2,3,4$, since $s_{k}$ is a Hermitian operator (as the sum
of products of commuting Hermitian operators), we know that 
\begin{equation}
0\leq\left(\mathrm{Tr}\left(Ws_{k}\right)\right)^{2}\leq\mathrm{Tr}\left(Ws_{k}^{2}\right).
\end{equation}
It follows from (\ref{eq:traces}) then that
\begin{equation}
\left|\mathrm{Tr}\left(W\left[a,a'\right]\left[b,b'\right]\right)\right|\leq4
\end{equation}
and

\begin{equation}
\mathrm{Tr}\left(Ws_{k}^{2}\right)\leq8.
\end{equation}
But then 
\begin{equation}
R_{k}^{2}=\left(\mathrm{Tr}\left(Ws_{k}\right)\right)^{2}\leq8.
\end{equation}
This implies (\ref{eq:Landau}) and (\ref{eq:Cirelson}).

That the value $2\sqrt{2}$ in (\ref{eq:Landau}) can be attained
is easy to show using the EPR/B paradigm: if $\alpha_{1}=0$, $\alpha_{2}=\pi/2$,
$\beta_{1}=\pi/4$, $\beta_{2}=-\pi/4$, then
\begin{equation}
R_{1}=\cos\left(\alpha_{1}-\beta_{1}\right)+\cos\left(\alpha_{1}-\beta_{2}\right)+\cos\left(\alpha_{2}-\beta_{1}\right)-\cos\left(\alpha_{2}-\beta_{2}\right)=2\sqrt{2}.
\end{equation}

\begin{rem}
It is instructive to see that if the operators $a,a'$ (or $b,b'$)
commute, (\ref{eq:traces}) leads to $R_{k}^{2}\leq4$, which, in
view of (\ref{eq:p-Exp}), is equivalent to (\ref{eq:Fine}). It is
tempting therefore to consider (\ref{eq:Fine}) as merely a special
(commutative) case of the construction used above to prove (\ref{eq:Cirelson}).
Notice however that this view cannot be accepted without additional
arguments: the proof of (\ref{eq:Fine}) makes no use of the assumption
that the outputs are eigenvalues of Hermitian operators in a Hilbert
space.
\end{rem}

\begin{rem}
It is known from {[}\ref{enu:Fine,-A.-(1982a).}, \ref{enu:Fine,-A.-(1982b).}{]}
that if a vector $\left(p_{11},p_{12},p_{21},p_{22}\right)$ satisfies
(\ref{eq:Fine}), then this vector can be generated by a system with
binary inputs and equiprobable binary outputs that satisfies (\ref{eq:selective influences}),
that is, is explainable by classical (non)contextuality. By contrast,
if a vector $\left(p_{11},p_{12},p_{21},p_{22}\right)$ satisfies
(\ref{eq:Cirelson}), it is not known to us whether this vector can
be generated by a quantum mechanical system with binary inputs and
equiprobable binary outputs. In this sense our characterization of
quantum contextuality is improvable. The issue of conditions that
are both necessary and sufficient for quantum contextuality has been
addressed {[}\ref{enu:Cabello-A-(2005)}, \ref{enu:Landau-LJ-(1988)},
\ref{enu:Masanes-Ll.-(2003)}{]}, but only in terms of the existence
of \emph{some} quantum systems, not necessarily those with binary
inputs and outputs. 
\end{rem}
\renewcommand{\thesection}{C}

\renewcommand{\thesubsection}{C\arabic{subsection}}

\setcounter{subsection}{0}

\setcounter{equation}{0}

\setcounter{thm}{0}

\section*{\label{sec:APPENDIX-C:-Computational}APPENDIX C: Computational Details}

\subsection{\label{sub:Computations-for-ELFP}Computations for ELFP}

A convex bounded polytope can be equivalently defined either as the
convex hull of a set of points (V-representation) or as the intersection
of half-spaces (H-representation). For our purposes, a V-representation
of a convex polytope in $d\textnormal{-}$space is given by a set
of points $x_{1},\dots,x_{n}\in\mathbb{R}^{d}$. The polytope consists
of all convex combinations of these points: $\lambda_{1}x_{1}+\dots+\lambda_{n}x_{n}$,
for all $\lambda_{1},\dots,\lambda_{n}\ge0$, $\lambda_{1}+\dots+\lambda_{n}=1$.
It is possible that the polytope is of lower dimension than the space
$\mathbb{R}^{d}$ in which it is defined if all the points $x_{i}$
reside in a lower dimensional affine subspace of $\mathbb{R}^{d}$.
A minimal V-representation (including only extreme points, i.e., points
that are vertices of the polytope) is unique. The H-representation
of a convex polytope is given by vectors $a_{1},\dots,a_{m}\in\mathbb{R}^{d}$
and a vector $b\in\mathbb{R}^{m}$. The polytope consists of the points
$x\in\mathbb{R}^{d}$ satisfying $a_{i}^{T}x\le b_{i}$ for all $i=1,\dots,m$.
A lower-dimensional convex polytope can be represented by including
inequalities of the forms $a^{T}x\le b$ and $(-a)^{T}x\le-b$ for
some $a$ and $b$ or by explicitly specifying certain constraints
as equations in the representation. For a full-dimensional convex
polytope, the minimal H-representation is unique. However, for a lower-dimensional
polytope, the equation constraints can be specified in many equivalent
ways and the inequality constraints can look different depending on
which of the linearly related coordinates are used to specify them.

There exist algorithms for converting between the two representations
of a convex polytope in exact rational arithmetic. We have used our
own program for these conversions but other programs, such as \emph{lrs}
(http://cgm.cs.mcgill.ca/\textasciitilde{}avis/C/lrs.html), can do
the same. The conversion between the two representation is computationally
demanding, the algorithms generally requiring superpolynomial time
in the size of the input.

A computationally simpler problem is eliminating redundant points
(those that are not vertices of the polytope) from a V-representation
or eliminating redundant equations or inequalities from an H-representation.
This problem can be solved by linear programming and the algorithm
is implemented in the \emph{redund} program that comes with \emph{lrs}.
However, the \emph{redund} program is not sufficient for putting an
H-representation to a minimal form as it cannot convert sets of inequalities
into equivalent equations (e.g., the three inequalities $x\ge0$,
$y\ge0$, $x+y\le0$ should be minimally represented as the two equations
$x=0$, $y=0$). To find the minimal H-representation, for every constraint
$a_{i}^{T}x\le b_{i}$ or $a_{i}^{T}x=b_{i}$ in turn, one can find
the upper and lower bounds $u$ and $l$ by maximizing and minimizing
the expression $a_{i}^{T}x$ given the other constraints, and apply
the following rules:
\begin{enumerate}
\item if this is an equation constraint (i.e., $a_{i}^{T}x=b_{i}$) and
$u=l=b_{i}$, then the constraint is redundant and can be eliminated; 
\item if this is an inequality constraint (i.e., $a_{i}^{T}x\le b_{i}$)
and $u\le b_{i}$, then the inequality is redundant and can be eliminated.
Otherwise, if $l=b_{i}$, then the constraint should be converted
to an equation.
\end{enumerate}
The dimension of a polytope can be determined from a minimal H-representation.
It is the dimension of the space minus the number of equation constraints
in the minimal representation. Given a full-dimensional polytope,
its volume can be computed using the \emph{lrs} program alongside
the conversion from a\textbf{ }V-representation to an H-representation.
If the polytope is given as an H-representation, then it has to be
converted to a V-representation first to compute its volume using
\emph{lrs}. To compute the volume of a lower-dimensional polytope,
we first move to a lower-dimensional parameterization that spans the
affine subspace where the polytope resides.

To compute ELFP, we begin by formulating the linear programming problem
$MQ=P$ subject to $Q\ge0$, as described in the main text ($M$ being
$2^{5}\times2^{8}$, $P$ having 2$^{5}$ components). $M$ defines
the V-representation for ELFP, and Vol$^{8}$ for ELFP is computed
directly from it. Applying an algorithm to find an equivalent H-representation
we obtain a system of 160 inequalities and 16 equations. We can then
substitute the expressions in the above matrices into this system
and reduce any redundant inequalities and equations. The resulting
system has 144 nonredundant inequalities and no equations with the
$p_{11},p_{12},p_{21},p_{22},\varepsilon_{1}^{1},\varepsilon_{2}^{1},\varepsilon_{1}^{2},\varepsilon_{2}^{2}$
variables. Then we algebraically simplify the list of 144 inequalities,
first into
\begin{equation}
\begin{array}{cc}
-\Gamma\le-p_{11}+p_{21}+p_{12}+p_{22} & \le1+\Gamma,\\
-\Gamma\le p_{11}-p_{21}+p_{12}+p_{22} & \le1+\Gamma,\\
-\Gamma\le p_{11}+p_{21}-p_{12}+p_{22} & \le1+\Gamma,\\
-\Gamma\le p_{11}+p_{21}+p_{12}-p_{22} & \le1+\Gamma,
\end{array}
\end{equation}

\begin{equation}
-\Lambda\le p_{11}+p_{21}+p_{12}+p_{22}\le2+\Lambda,
\end{equation}
 where
\begin{equation}
\begin{array}{cc}
|-p_{11}-p_{21}+p_{12}+p_{22}| & \le1+\Lambda,\\
|-p_{11}+p_{21}-p_{12}+p_{22}| & \le1+\Lambda,\\
|-p_{11}+p_{21}+p_{12}-p_{22}| & \le1+\Lambda,
\end{array}
\end{equation}
\begin{equation}
\begin{array}{cc}
\Gamma=\min\{\, & 1-\varepsilon_{1}^{1}-\varepsilon_{1}^{2}+\varepsilon_{2}^{1}+\varepsilon_{2}^{2},\\
 & 1-\varepsilon_{1}^{1}+\varepsilon_{1}^{2}-\varepsilon_{2}^{1}+\varepsilon_{2}^{2},\\
 & 1-\varepsilon_{1}^{1}+\varepsilon_{1}^{2}+\varepsilon_{2}^{1}-\varepsilon_{2}^{2},\\
 & 1+\varepsilon_{1}^{1}-\varepsilon_{1}^{2}-\varepsilon_{2}^{1}+\varepsilon_{2}^{2},\\
 & 1+\varepsilon_{1}^{1}-\varepsilon_{1}^{2}+\varepsilon_{2}^{1}-\varepsilon_{2}^{2},\\
 & 1+\varepsilon_{1}^{1}+\varepsilon_{1}^{2}-\varepsilon_{2}^{1}-\varepsilon_{2}^{2},\\
 & \varepsilon_{1}^{1}+\varepsilon_{1}^{2}+\varepsilon_{2}^{1}+\varepsilon_{2}^{2},\\
 & 2-\varepsilon_{1}^{1}-\varepsilon_{1}^{2}-\varepsilon_{2}^{1}-\varepsilon_{2}^{2}\,\},
\end{array}
\end{equation}
\begin{equation}
\begin{array}{cc}
\Lambda=\min\{\,- & \varepsilon_{1}^{1}+\varepsilon_{1}^{2}+\varepsilon_{2}^{1}+\varepsilon_{2}^{2},\\
 & \varepsilon_{1}^{1}-\varepsilon_{1}^{2}+\varepsilon_{2}^{1}+\varepsilon_{2}^{2},\\
 & \varepsilon_{1}^{1}+\varepsilon_{1}^{2}-\varepsilon_{2}^{1}+\varepsilon_{2}^{2},\\
 & \varepsilon_{1}^{1}+\varepsilon_{1}^{2}+\varepsilon_{2}^{1}-\varepsilon_{2}^{2},\\
 & 1-\varepsilon_{1}^{1}-\varepsilon_{1}^{2}-\varepsilon_{2}^{1}+\varepsilon_{2}^{2},\\
 & 1-\varepsilon_{1}^{1}-\varepsilon_{1}^{2}+\varepsilon_{2}^{1}-\varepsilon_{2}^{2},\\
 & 1-\varepsilon_{1}^{1}+\varepsilon_{1}^{2}-\varepsilon_{2}^{1}-\varepsilon_{2}^{2},\\
 & 1+\varepsilon_{1}^{1}-\varepsilon_{1}^{2}-\varepsilon_{2}^{1}-\varepsilon_{2}^{2}\,\},
\end{array}
\end{equation}
 and then, by noticing regularities, into the compact inequality (\ref{eq:ELFT}). 
\begin{rem}
Changing $\varepsilon_{j}^{i}\rightarrow\nicefrac{1}{2}-\varepsilon_{j}^{i}$
leads to (denoting the new $\varepsilon\textnormal{-}$vector by $\varepsilon'$)
\begin{equation}
\max\mathsf{S}_{1}\varepsilon'=\max\mathsf{S}_{0}\varepsilon,\max\mathsf{S}_{0}\varepsilon'=\max\mathsf{S}_{1}\varepsilon.
\end{equation}
Analogously for $p_{ij}\rightarrow\nicefrac{1}{2}-p_{ij}$, 
\begin{equation}
\max\mathsf{S}_{1}p'=\max\mathsf{S}_{0}p,\max\mathsf{S}_{0}p'=\max\mathsf{S}_{1}p.
\end{equation}
It follows that we cannot without loss of generality confine all components
of $\varepsilon$ or $p$ to $\left[0,\nicefrac{1}{4}\right]$. But
ELFP does not change if the transformation $x\rightarrow\nicefrac{1}{2}-x$
is applied to an even number of the components of $\left(p,\varepsilon\right)$. 
\end{rem}

\subsection{\label{sub:Computations-for-3,}Computations for $\mathrm{chaos(p)}$,
$\mathrm{quant}(p)$, and $\mathrm{class}(p)$ constraints}

The All$_{\mathrm{constr}}$ polytopes for the three constraints are
obtained by concatenating the ELFP equations and inequalities with
the constraint inequalities. Then, the volumes are computed by using
the \emph{lrs} program as described above.

For Fit$_{\mathrm{constr}}$ polytopes, we observe first that they
are convex. This follows from
\begin{equation}
\begin{array}{r}
\mathrm{Fit_{constr}}=\{\varepsilon:\forall i=1,\dots,n:\left(p_{\left(i\right)},\varepsilon\right)\in\mathrm{ELFP\,}\}\\
=\mathrm{ELFP}_{p_{\left(1\right)}}\cap\cdots\cap\mathrm{ELFP}_{p_{\left(n\right)}},
\end{array}
\end{equation}
where $p_{\left(i\right)}$, $i=1,\dots,n$, denote the vertices of
the 4D convex polytope defined by $\mathrm{constr}$ and $\mathrm{ELFP}_{p_{\left(i\right)}}$
denotes the (convex) cross-section of the ELFP set formed with $p=p_{\left(i\right)}$.
It follows that $\mathrm{Fit_{\mathrm{constr}}}$ is convex as the
intersection of convex sets. Following the logic of this observation,
we have implemented a general program for eliminating variables from
a system of linear equations and inequalities so that the resulting
system is satisfied for exactly those values for which there exist
such values of the eliminated variables for which the original system
is satisfied. This program together with steps to ensure that the
resulting representation is minimal was used to find all the Fit sets
shown in the main text.

Finding the forcing sets is more difficult as they are generally not
convex. We characterize them using the equation 
\begin{equation}
\begin{array}{l}
\mathrm{Force}_{\mathrm{chaos}}-\mathrm{Force}_{\mathrm{constr}}\\
=\left\{ \varepsilon:(\exists p:\left(p,\varepsilon\right)\in\mathrm{ELFP}\wedge\neg\mathrm{constr}\left(p\right))\right\} .
\end{array}
\end{equation}
This equation provides an algorithm: for each inequality in $\mathrm{constr}$,
form the conjunction of the ELFP inequalities with the negation of
the inequality. Then project this conjunction to the $\varepsilon$
4-space. The union of these projections over all inequalities in $\mathrm{constr}$
is the set $\mathrm{Force}_{\mathrm{chaos}}-\mathrm{Force}_{\mathrm{constr}}$.
We have implemented a general program that takes as input a representation
of a polytope, a list of additional constraints, and a list of variables
to eliminate. It then outputs a representation of the difference of
the polytope and the set represented by the additional constraints
projected to the remaining (not eliminated) variables. This representation
consists of a list of linear systems whose disjunction characterizes
the resulting set. In all our computations it turned out that all
the linear systems in the disjunction were the same, and so the sets
$\mathrm{Force}_{\mathrm{chaos}}-\mathrm{Force}_{\mathrm{constr}}$
are in fact convex in these cases.

The computations of Equi sets require no elaboration.
\begin{rem}
There is the practical problem that the negation of a $\le\textnormal{-}$inequality
is a $>\textnormal{-}$inequality while standard algorithms only accept
closed convex polytopes. To cope with this problem, we approximated
$a>b$ by $a\ge b+\text{(very small number)}$. We also used a rational
approximation to $\sqrt{2}$ in the quant constraints. In both cases,
we have repeated the computations with decreasing values of ``very
small number'' until it was obvious where the results converged.
\end{rem}

\subsection{\label{sub:Computations-for-1}Computations for $\mathrm{Fit}_{\mathrm{fix}\left(p\right)}$
constraint}

That $\max\mathsf{S}_{0}p$ and $\max\mathsf{S}_{1}p$ are contained
in and completely fill the triangle $\left\{ (0,0),(\nicefrac{1}{2},1),(1,\nicefrac{1}{2})\right\} $
can be verified by splitting\textbf{ }(\ref{eq:fixed}) into 64 component
cases according as which of the values of $\mathsf{S}_{0}p$ and $\mathsf{S}_{1}p$
are the maxima, finding the vertices of each component system, and
drawing the union of these components in $\max\mathsf{S}_{0}p$ and
$\max\mathsf{S}_{1}p$ coordinates. The triangle is described by
\begin{equation}
\begin{array}{c}
2\max\mathsf{S}_{0}p-\max\mathsf{S}_{1}p\ge0,\\
2\max\mathsf{S}_{1}p-\max\mathsf{S}_{0}p\ge0,\\
\max\mathsf{S}_{0}p+\max\mathsf{S}_{1}p\le\nicefrac{3}{2}.
\end{array}
\end{equation}
 Adding these inequalities to the representation of (\ref{eq:fixed})
as linear inequalities according to the definitions of $\max\mathsf{S}_{0}\varepsilon$
and $\max\mathsf{S}_{1}\varepsilon$, we obtain a 6D polytope $P^{\left(6\right)}$
in $\left(\varepsilon,\max\mathsf{S}_{0}p,\max\mathsf{S}_{1}p\right)\textnormal{-}$coordinates.\textbf{
}In the V-representation of $P^{\left(6\right)}$, all vertices have
values of $\max\mathsf{S}_{0}p$ and $\max\mathsf{S}_{1}p$ in the
set 
\begin{equation}
\left\{ (0,0),(\nicefrac{1}{4},\nicefrac{1}{2}),(\nicefrac{1}{2},\nicefrac{1}{4}),(\nicefrac{1}{2},1),(1,\nicefrac{1}{2})\right\} .
\end{equation}
It follows that every edge of the polytope projects to one of these
5 points or to a line connecting two of them. Consequently, as $\left(\max\mathsf{S}_{0}\varepsilon,\max\mathsf{S}_{1}\varepsilon\right)$
changes within any triangle $T$ formed by these lines, the cross-section
$P_{\left(\max\mathsf{S}_{0}\varepsilon,\max\mathsf{S}_{1}\varepsilon\right)}^{\left(4\right)}$
of $P^{\left(6\right)}$ retains its structure (face lattice) while
its coordinates change as affine functions of $\left(\max\mathsf{S}_{0}\varepsilon,\max\mathsf{S}_{1}\varepsilon\right)\in T$.
It follows that the volume of $P_{\left(\max\mathsf{S}_{0}\varepsilon,\max\mathsf{S}_{1}\varepsilon\right)}^{\left(4\right)}$
is a polynomial of $\left(\max\mathsf{S}_{0}\varepsilon,\max\mathsf{S}_{1}\varepsilon\right)\in T$
of at most degree four.\textbf{ }The coefficients of these polynomials
were obtained by fitting unconstrained degree 4 polynomials to the
exact volumes $\mathrm{Vol}^{4}\left(\mathrm{Fit}_{\mathrm{fix}\left(p\right)}\right)$
for $\left(\max\mathsf{S}_{0}p,\max\mathsf{S}_{1}p\right)\in\left\{ 0,.01,.02,\dots,\ensuremath{1}\right\} ^{2}$.\textbf{
}It turns out that the coefficients change only if either of the differences
$\max\mathsf{S}_{0}p-\nicefrac{1}{2}$ and $\max\mathsf{S}_{1}p-\nicefrac{1}{2}$
changes its sign. In all cases the fit is perfect for the number of
points far exceeding the number of coefficients, confirming that the
computations are correct.

\pagebreak{}

\section*{REFERENCES}
\begin{enumerate}
\item \label{enu:Townsend,-J.T.,-&}Townsend JT, Schweickert R (1989) Toward
the trichotomy method of reaction times: Laying the foundation of
stochastic mental networks.\emph{ J. Math. Psych}. \textbf{33}: 309-327.
\item \label{enu:Bell,-J.-(1964).}Bell J (1964) On the Einstein-Podolsky-Rosen
paradox. \emph{Physics} \textbf{1}: 195-200.
\item \label{enu:J.-Bell,-1966}Bell J (1966) On the problem of hidden variables
in quantum mechanics. \emph{Rev. Modern Phys.} \textbf{38}: 447-453.
\item \label{enu:Clauser,-J.-F.,1969}Clauser JF, Horne MA, Shimony A, Holt
RA (1969) Proposed experiment to test local hidden-variable theories.
\emph{Phys. Rev. Lett.} \textbf{23}: 880-884.
\item \label{enu:Clauser,-J.F.-and}Clauser JF, Horne MA (1974) Experimental
consequences of objective local theories\emph{. Phys. Rev. D}, \textbf{10}:
526-535.
\item \label{enu:Fine,-A.-(1982a).}Fine A (1982) Joint distributions, quantum
correlations, and commuting observables.\emph{ J. Math. Phys}. \textbf{23}:
1306-1310.
\item \label{enu:Fine,-A.-(1982b).}Fine A (1982) Hidden variables, joint
probability, and the Bell inequalities.\emph{ Phys. Rev. Lett.} \textbf{48}:
291-295.
\item \label{enu:Dzhafarov,-E.N.-(2003).}Dzhafarov EN (2003) Selective
influence through conditional independence.\emph{ Psychom}. \textbf{68}:
7-26.
\item \label{enu:D.-Bohm,-&}Bohm D, Aharonov Y (1957) Discussion of experimental
proof for the paradox of Einstein, Rosen and Podolski. \emph{Phys.
Rev}., \textbf{108}: 1070-1076.
\item \label{enu:E.N.-Dzhafarov,-&2012}Dzhafarov EN, Kujala JV (2012) Selectivity
in probabilistic causality: Where psychology runs into quantum physics.
\emph{J. Math. Psych}., \textbf{56}: 54-63.
\item \label{enu:E.N.-Dzhafarov,-in press 2}Dzhafarov EN, Kujala JV (2012)
Quantum entanglement and the issue of selective influences in psychology:
An overview. \emph{Lect. Notes in Comp. Sci}. 7620, 184-195.
\item \label{enu:J.-Cereceda,-Quantum}Cereceda J (2000) Quantum mechanical
probabilities and general probabilistic constraints for Einstein\textendash{}Podolsky\textendash{}Rosen\textendash{}Bohm
experiments. \emph{Found. Phys. Lett.} \textbf{13}: 427-442
\item \label{enu:Kujala,-J.-V.,2008}Kujala JV, Dzhafarov EN (2008) Testing
for selectivity in the dependence of random variables on external
factors.\emph{ J. Math. Psych}. \textbf{52}: 128-144.
\item \label{enu:Dzhafarov,-E.N.,-&2010}Dzhafarov EN, Kujala JV (2010)
The Joint Distribution Criterion and the Distance Tests for selective
probabilistic causality.\emph{ Front. Quant. Psych. Meas}. \textbf{1:151}
\href{http://www.frontiersin.org/Quantitative_Psychology_and_Measurement/10.3389/fpsyg.2010.00151/full}{doi: 10.3389/fpsyg.2010.00151}.
\item \label{enu:E.N.-Dzhafarov,-&in press}Dzhafarov EN, Kujala JV (in
press) Order-distance and other metric-like functions on jointly distributed
random variables. \emph{Proc. Amer. Math. Soc}. (available as \href{http://arxiv.org/pdf/1110.1228v3.pdf}{arXiv:1110.1228v3 [math.PR]})
\item \label{enu:E.N.-Dzhafarov,-R.Sung 2004}Dzhafarov EN, Schweickert
R, Sung K (2004) Mental architectures with selectively influenced
but stochastically interdependent components. \emph{J. Math. Psych}.,
\textbf{48}: 51-64.
\item \label{enu:Schweickert,-R.,-Fisher,}Schweickert R, Fisher DL, Goldstein
WM (2010) Additive factors and stages of mental processes in task
networks.\emph{ J. Math. Psych.} \textbf{54}: 405-414.
\item \label{enu:S.-Sternberg,-The}Sternberg S (1969) The discovery of
processing stages: Extensions of Donders\textquoteright{} method.
In W.G. Koster (Ed.), \emph{Attention and Performance II. Acta Psych}.,
\textbf{30}: 276\textendash{}315.
\item \label{enu:J.-T.-Townsend,1984}Townsend JT (1984) Uncovering mental
processes with factorial experiments. \emph{J. Math. Psych}., \textbf{28}:
363\textendash{}400.
\item \label{enu:P.-Suppes,-M.1981}Suppes P, Zanotti M (1981) When are
probabilistic explanations possible?\emph{ Synthese} \textbf{48}:
191-199. 
\item \label{enu:P.-Suppes,-Representation}Suppes P (2002) Representation\emph{
and Invariance of Scientific Structures} (CSLI, Stanford) pp. 332-351.
\item \label{enu:A.-Aspect,-P.1981}Aspect A, Grangier P, Roger G (1981)
Experimental tests of realistic local theories via Bell's theorem.
\emph{Phys. Rev. Lett.} \textbf{47}: 460\textendash{}463.
\item \label{enu:A.-Aspect,-P.1982}Aspect A, Grangier P, Roger G (1982)
Experimental realization of Einstein-Podolsky-Rosen-Bohm gedankenexperiment:
A new violation of Bell's inequalities. \emph{Phys. Rev. Lett.} \textbf{49}:
91\textendash{}94.
\item \label{enu:Cirel'son-BS-(1980)}Cirel'son BS (1980) Quantum generalizations
of Bell's inequality.\emph{ Lett. Math. Phys.} \textbf{4}: 93--100.
\item \label{enu:Landau-LJ-(1987)}Landau LJ (1987) On the violation of
bell's inequality in quantum theory. \emph{Phys. Lett. A} \textbf{120:}
54--56.
\item \label{enu:H.-Thorisson,-Coupling,}Thorisson H (2000) \emph{Coupling,
Stationarity, and Regeneration} (New York, Springer).
\item \label{enu:F.-Laudisa,-Contextualism}Laudisa F (1997) Contextualism
and nonlocality in the algebra of EPR observables. \emph{Phil. Sci.},
\textbf{64}: 478-496. 
\item \label{enu:A.-Yu.-Khrennikov,}Khrennikov AYu (2008) EPR\textendash{}Bohm
experiment and Bell\textquoteright{}s inequality: Quantum physics
meets probability theory. \emph{Theor. Math. Phys.}, \textbf{157}:
1448\textendash{}1460.
\item \label{enu:A.-Yu.-Khrennikov,-1}Khrennikov AYu (2009) \emph{Contextual
Approach to Quantum Formalism. Fundamental Theories of Physics} \textbf{160}
(Dordrecht, Springer).
\item \label{enu:I.-Pitowsky,-Quantum}Pitowsky I \emph{(}1989)\emph{ Quantum
Probability--Quantum Logic, Lect. Notes in Phys.} \textbf{321} (Heidelberg,
Springer). 
\item \label{enu:A.-Peres,-All}Peres A (1999) All the Bell inequalities.
\emph{Found. Phys.} \textbf{29} :589-614.
\item \label{enu:R.M.-Basoalto,-&}Basoalto RM, Percival IC (2003) BellTest
and CHSH experiments with more than two settings. \emph{J. Phys. A:
Math. \& Gen.} \textbf{36}: 7411\textendash{}7423.
\item \label{enu:Cabello-A-(2005)}Cabello A (2005) How much larger quantum
correlations are than classical ones. Phys. Rev. A \textbf{72}: 12113
1-5.
\item \label{enu:Landau-LJ-(1988)}Landau LJ (1988) Empirical two-point
correlation functions. \emph{Found. Phys.} \textbf{18}: 449--460.
\item \label{enu:Masanes-Ll.-(2003)}Masanes Ll. (2003) Necessary and sufficient
condition for quantum-generated correlations. \emph{arXiv:quant-ph/0309137}.\end{enumerate}

\end{document}